\documentclass[11pt, a4paper]{article}

\usepackage{amssymb,amsmath, wasysym, graphicx, epsfig, setspace, wrapfig, comment}
\usepackage{latexsym,dsfont,amsfonts,amsmath,amssymb,amsthm}
\usepackage{pstricks,epsfig,psfrag,color}
\usepackage{bbm}
\usepackage[show]{ed}
\usepackage{soul}
\usepackage{mathrsfs}   
\usepackage{tikz}
\usepackage{upgreek}

\usepackage{mathtools}
\usepackage{graphicx}
\usepackage{bbm}
\usepackage{subfig} 
\usepackage{amsmath,amsfonts,amssymb}  

\numberwithin{equation}{section}

\usepackage{bm}

\usepackage{float}
\usepackage{mathtools}
\usepackage{hyperref}
\usepackage[capitalize]{cleveref}
\crefname{equation}{}{}

\usepackage[style=numeric,
            sorting=nyt,
            giveninits=true,
            maxbibnames=20,
            natbib=true,
            backend=bibtex]{biblatex}
\AtEveryBibitem{%
  \clearfield{issn} 
  \clearfield{month} 
  \clearfield{language} 

  \ifentrytype{online}{}{
    \clearfield{url}
  }
}
\usepackage{xurl}	

\usepackage{color}
\usepackage{comment}
\usepackage{tikz-cd}
\usepackage[shortlabels]{enumitem}
\usepackage{nomencl}

\newtheorem{theorem}{Theorem}[section]
\newtheorem{lemma}[theorem]{Lemma}
\newtheorem{corollary}[theorem]{Corollary}

\theoremstyle{definition}
\newtheorem{remark}[theorem]{Remark}
\theoremstyle{plain}

\numberwithin{equation}{section}

\usepackage{amsfonts}
\usepackage{graphicx}
\usepackage{epstopdf}
\usepackage{algorithmic}
\ifpdf
  \DeclareGraphicsExtensions{.eps,.pdf,.png,.jpg}
\else
  \DeclareGraphicsExtensions{.eps}
\fi

\title{Robust spectral preconditioning for high-P\'{e}clet number convection-diffusion}

\author{Lukas Holbach\footnotemark[1]
\and Peter Bastian\footnotemark[1]
\and Robert Scheichl\footnotemark[2] $^{,}$\footnotemark[1]}

\usepackage{amsopn}

\newcommand{\N}{{\mathord{\mathbb N}}}
\newcommand{\R}{{\mathord{\mathbb R}}}

\def\DHLhksqrt#1#2{%
\setbox0=\hbox{$#1\sqrt{#2\,}$}\dimen0=\ht0
\advance\dimen0-0.2\ht0
\setbox2=\hbox{\vrule height\ht0 depth -\dimen0}%
{\box0\lower0.4pt\box2}}

\renewcommand{\vec}[1]{\mathbf{#1}}

\newcommand{\dx}{\ensuremath{\, d\vec{x}}}
\newcommand{\ds}{\ensuremath{\, ds}}

\makeatletter
\providecommand*{\cupdot}{%
  \mathbin{%
    \mathpalette\@cupdot{}%
  }%
}
\newcommand*{\@cupdot}[2]{%
  \ooalign{%
    $\m@th#1\cup$\cr
    \hidewidth$\m@th#1\cdot$\hidewidth
  }%
}
\makeatother

\makeatletter
\newcommand{\subalign}[1]{%
  \vcenter{%
    \Let@ \restore@math@cr \default@tag
    \baselineskip\fontdimen10 \scriptfont\tw@
    \advance\baselineskip\fontdimen12 \scriptfont\tw@
    \lineskip\thr@@\fontdimen8 \scriptfont\thr@@
    \lineskiplimit\lineskip
    \ialign{\hfil$\m@th\scriptstyle##$&$\m@th\scriptstyle{}##$\crcr
      #1\crcr
    }%
  }
}
\makeatother

\newcommand{\nDoms}{M}
\newcommand{\nOvs}{\ell}
\newcommand{\ddindex}{j}
\newcommand{\locDim}{n}

\newcommand{\ovlpdomain}{\omega}
\newcommand{\colConst}{\kappa}
\newcommand{\osdomain}{\ovlpdomain^*}
\newcommand{\osbdry}{\partial \osdomain_\ddindex}
\newcommand{\pu}{\chi}
\newcommand{\puC}{C_\chi}
\newcommand{\resOp}{P}
\newcommand{\eval}{\lambda}
\newcommand{\evec}{\varphi}
\newcommand{\evalMax}{\eval_{\mathrm{max}}}

\newcommand{\evalFixed}{n_{\mathrm{sd}}}
\newcommand{\evalindex}{k}
\newcommand{\evalloc}[1][\evalindex]{\eval_{\ddindex,#1}}
\newcommand{\evecloc}[1][\evalindex]{\evec_{\ddindex,#1}}
\newcommand{\nIter}{\#IT}
\newcommand{\errTol}{\varepsilon}

\newcommand{\fwpTestD}[1]{H^1_{D}(#1)}
\newcommand{\fwpTestSpace}{\fwpTestD{\domain}}
\newcommand{\fwpTestSpaceLoc}[1][\osdomain]{H^1_{D,0}(#1_\ddindex)}
\newcommand{\harmSpace}{H_{B}(\osdomain_\ddindex)}
\newcommand{\locApproxSpace}{S_{\locDim_\ddindex}(\ovlpdomain_\ddindex)}
\newcommand{\globalApproxSpace}{S_{\locDim}(\domain)}
\newcommand{\msgfemOperator}{G}
\newcommand{\projB}{\pi_\ddindex}
\newcommand{\projCoarse}{\pi_S}
\newcommand{\discRestriction}{\mat{R}_\ddindex}
\newcommand{\discCoarseRestriction}{\mat{R}_S}
\newcommand{\discBilForm}{\mat{B}}
\newcommand{\discPU}{\bs{\pu}_\ddindex}
\newcommand{\discProjB}{\bs{\pi}_\ddindex}
\newcommand{\discProjCoarse}{\bs{\pi}_S}
\newcommand{\discPreconditioner}{\mat{M}}

\newcommand{\face}{\gamma}
\newcommand{\element}{\tau}
\newcommand{\mesh}{\mathcal{T}_h}

\newcommand{\diffTensor}{A}
\newcommand{\diffMin}{a_{\mathrm{min}}}
\newcommand{\diffMax}{a_{\mathrm{max}}}
\newcommand{\velField}{\vec{b}}
\newcommand{\fwp}{u}
\newcommand{\fwpTest}{v}
\newcommand{\fwpSol}{\fwp^e}
\newcommand{\fwpGFEM}{\fwp^G}
\newcommand{\fwpPart}{\fwp^p}
\newcommand{\fwpDBC}{\fwp_\bcD}
\newcommand{\fwpiPart}{\fwp_\ddindex^p}
\newcommand{\rhssource}{f}
\newcommand{\bcD}{g}
\newcommand{\bcO}{q}
\newcommand{\peclet}{\mathrm{P\acute{e}}}
\newcommand{\bilB}{B}
\newcommand{\bilWF}[3]{\bilB_{#3}\mspace{-2mu}\left(#1,#2\right)}
\newcommand{\bilIP}[3]{a_{#3}\mspace{-2mu}\left(#1,#2\right)}
\newcommand{\rhsF}{F}
\newcommand{\rhsWF}[2]{\rhsF_{#2}\mspace{-2mu}\left(#1\right)}

\newcommand{\domain}{\Omega}
\newcommand{\boundary}{\partial \domain}
\newcommand{\bdryD}{\Gamma_D}
\newcommand{\bdryO}{\Gamma_O}

\newcommand{\x}{\vec{x}}

\newcommand{\Rd}{\R^\dimDomain}

\newcommand{\Lp}[1]{L^{#1}(\domain)}
\newcommand{\Lpd}[1]{L^{#1}(\domain;\,\R^\dimDomain)}

\newcommand{\Lpdd}[1]{L^{#1}(\domain;\,\R^{\dimDomain \times \dimDomain})}

\newcommand{\sobolev}[2]{W^{#1,#2}(\domain)}
\newcommand{\sobolevD}[2]{W^{#1,#2}(\domain;\R^\dimDomain)}

\newcommand{\Hone}{H^1(\domain)}

\newcommand{\dimDomain}{d}

\newcommand{\trans}{\mathsf{T}}
\newcommand{\transpose}{^{\trans}}

\newcommand{\dual}{^{\prime}}

\DeclareMathOperator*{\divergence}{div}

\newcommand{\abs}[1]{\left\lvert#1\right\rvert}

\newcommand{\norm}[1]{\left\lVert#1\right\rVert}

\DeclareMathOperator{\diam}{diam}

\DeclareMathOperator{\supp}{supp}

\newcommand{\bs}[1]{{\boldsymbol #1}}
\newcommand{\boldvar}[1]{\ensuremath{\boldsymbol{#1}}}
\newcommand{\boldvec}[1]{\ensuremath{\mathbf{#1}}}
\newcommand{\mat}[1]{\boldvec{#1}}
\ifdefined\vec
  \renewcommand{\vec}[1]{\boldvec{#1}}
\else
  \newcommand{\vec}[1]{\boldvec{#1}}
\fi

\newcommand{\bn}{\boldvar{\nu}}

\usepackage{booktabs}
\usepackage{subcaption}
\usepackage{stmaryrd}
\usepackage{tikz}
\usepackage[shortlabels]{enumitem}

\date{\today}

\makeatletter
\let\@fnsymbol\@arabic
\makeatother

\begin{document}

\maketitle

\footnotetext[1]{Interdisciplinary Center for Scientific Computing, Heidelberg University, Germany (\texttt{lukas.holbach@iwr.uni-heidelberg.de}, \texttt{peter.bastian@iwr.uni-heidelberg.de}).}
\footnotetext[2]{Institute for Mathematics, Heidelberg University, Germany (\texttt{r.scheichl@uni-heidelberg.de}).}

\begin{abstract}
We introduce a two-level hybrid restricted additive Schwarz (RAS) preconditioner for heterogeneous steady-state convection-diffusion equations at high P\'{e}clet numbers. Our construction builds on the multiscale spectral generalized finite element method (MS-GFEM), wherein the coarse space is spanned by locally optimal basis functions obtained from local generalized eigenproblems on operator-harmonic spaces. Extending the theory of \cite{Ma2025} to convection-diffusion problems in conservation form, we establish exponential convergence of the MS-GFEM approximation with respect to the dimension of the local approximation space. Rewriting MS-GFEM as a RAS‐type iteration, we show for coercive problems that this exponential convergence property is inherited by the RAS-type iterative method (at least in the continuous setting). Employed as a preconditioner within the generalized minimal residual method (GMRES), the resulting method requires only a few iterations for high accuracy even with low‐dimensional coarse spaces.

Through extensive numerical experiments on problems with high‐contrast diffusion and non-divergence-free, rotating velocity fields, we demonstrate robustness with respect to the grid P\'{e}clet number and the number of subdomains (tested up to $10^5$ subdomains), while coarse‐space dimensions remain small as grid P\'{e}clet numbers increase. By adapting the coarse space and oversampling size, we are able to achieve arbitrarily fast convergence of preconditioned GMRES. As an extension, for which we do not have theory yet, we show effectiveness of the method even for indefinite problems and in the vanishing‐diffusion limit.
\end{abstract}



\section{Introduction}

Steady-state convection-diffusion equations play an important role in the modeling of transport phenomena across a wide range of scientific and engineering disciplines. In groundwater hydrology, convection-diffusion equations form the basis for geostatistical inversion methods used to infer subsurface properties from tracer experiments \cite{CirpkaKitanidis2000,NgoBastianIppisch2015}. In semiconductor device simulation, they describe charge carrier transport under the influence of electric fields and diffusion processes \cite{PolakEtAl1987,ChenBagci2020}. Similar mathematical structures arise in the study of stellarator plasmas, where they model the steady transport of heat and particles in magnetically confined fusion devices \cite{Escoto2024}. 
These diverse applications highlight the relevance of developing efficient and robust numerical methods for solving stationary convection-diffusion problems, which also serve as building blocks for implicit time discretizations of unsteady systems with large time steps.

In many practical settings, flow is dominated by convection, i.e.,  the P\'{e}clet number is high. At the same time, the diffusion tensor may be highly heterogeneous and exhibit large contrasts, causing sharp internal or boundary layers in the solution. Standard conforming finite-element discretizations then require prohibitively fine meshes to resolve these layers and frequently produce oscillatory or inaccurate approximations \cite{Stynes2005}. While there are several approaches to alleviating these issues (see e.g.\ the text book \cite{RoosStynesTobiska2010} and \cite{Stynes2005}), even with stabilized discretizations such as streamline upwind Petrov-Galerkin (SUPG) or discontinuous Galerkin (DG) methods, the resulting linear systems are non-symmetric and typically very large. Hence, robust iterative solvers and preconditioners are essential. 

Among the most successful preconditioners for symmetric positive-definite (SPD) sparse linear systems are multigrid methods. Choosing carefully designed components they can also be applied to stationary convection-diffusion problems.
Using robust smoothers to achieve P\'{e}clet number robustness has been investigated in \cite{BEY1997177,HackbuschProbst1997,NgoBastianIppisch2015} where the almost upper triangular structure of the system matrix of upwind schemes is exploited in ``downwind numbering''. These schemes, however, are hard to parallelize efficiently on many processors \cite{Manteuffel2019}. 
Robust coarse grid corrections in the context of algebraic multigrid have been investigated recently in a series of papers \cite{Manteuffel2018,Manteuffel2019,Ali2024}.

The preconditioner we propose falls into the class of domain decomposition methods (DDMs), where the global problem is divided into smaller local subproblems that can be solved independently in parallel. Subsequently, a global solution is constructed by combining the local solutions, e.g., using a partition of unity, and this procedure is repeated iteratively with updated local boundary conditions. One-level DDMs are easy to implement, but iteration counts typically grow with the number of subdomains, since information propagates only between neighboring subdomains in each iteration (see e.g.\ \cite[Chapter 4]{DoleanJolivetNataf2015}). Introducing a global coarse-space correction yields two-level methods that can restore scalability, but the design of coarse spaces that are robust with respect to heterogeneous coefficients and other model parameters is delicate.
Classical coarse spaces are built as piecewise polynomials on a coarse grid \cite{CaiWidlund1992,CaiWidlund1993} but are not robust with respect to heterogeneities in model parameters. Other examples include the finite element tearing and interconnecting (FETI) method \cite{Toselli2001}, the balancing domain decomposition by constraints (BDDC) method \cite{TuLi2008}, as well as Neumann-Neumann or Robin-Robin methods \cite{Achdou2000}. For a more detailed overview of classical DDMs for convection-diffusion equations, we refer the reader to \cite[Section 11.5.1]{ToselliWidlund2005} and the references therein.

A common approach for obtaining parameter-robust coarse spaces is to build them from eigenfunctions of local generalized eigenvalue problems, which identify problem-dependent low-energy modes that must be represented in the coarse space. Such spectral coarse spaces were introduced in the context of BDDC and FETI methods in \cite{MANDEL20071389} and for overlapping Schwarz methods in \cite{GalvisEfendiev2010a,GalvisEfendiev2010b}. The generalized eigenvalue problems on the overlap (GenEO) coarse space, which was originally proposed for SPD systems in \cite{Spillane2013}, has been applied to non-self-adjoint and indefinite problems in \cite{Bootland2022} by using a generalized eigenvalue problem based on a self-adjoint and positive-definite subproblem. However, results were not fully robust for high P\'{e}clet numbers. The authors of \cite{AlDaasJolivetRees2023} propose a fully algebraic variant of GenEO based on complex generalized eigenproblems. They analyze the method for diagonally dominant SPD problems and demonstrate its performance for non-SPD problems in numerical examples. Another extension of GenEO for non-symmetric systems, which reduces to our method in the SPD case, is presented in \cite{NatafParolin2025}. A generalized Dryja-Smith-Widlund-type coarse space for incompressible fluid flow has been presented in \cite{gdsw_navierstokes_2020}.

Similar ideas of spectral approximation spaces for convection-diffusion equations have been employed in the context of multiscale PDEs, e.g., local orthogonal decomposition \cite{LiPeterseimSchedensack2017,Bonizzoni2024} or multiscale finite element methods \cite{Calo2016,FuChungLi2025}. While the goal of multiscale methods is to incorporate important physical microscale features into a coarse macroscopic ansatz space for the solution of the PDE, these approximation spaces can still become prohibitively large when accurate solutions are needed. Hence, it is a natural idea to combine such coarse spaces with DDMs to build an iterative method that requires significantly smaller coarse spaces while maintaining accuracy in the solutions. 

The method we propose is based on the multiscale spectral generalized finite element method (MS-GFEM), which was originally introduced for elliptic problems in \cite{BabuskaLipton2011} and then further developed in \cite{Ma2022}. Recently, a unified framework for MS-GFEM that is applicable to a wide range of multiscale PDEs with $L^\infty$-coefficients has been presented in \cite{Ma2025}. The MS-GFEM approximation space is built from local generalized eigenvalue problems, where the eigenfunctions are sought in a problem-dependent local operator-harmonic subspace. With this construction, a nearly exponential decay of the eigenvalues can be proved, which is crucial for the efficiency of the method, since the global error bound depends on the eigenvalues. In this paper, we incorporate exactly this MS-GFEM approximation space as a coarse space into a two-level hybrid restricted additive Schwarz preconditioner. Exploiting the rapid eigenvalue decay allows us to achieve robust preconditioning with small coarse spaces. Note that similar ideas have recently been applied to elliptic \cite{Strehlow2024} and Helmholtz problems \cite{Ma2025Helmholtz}.

Our main contributions are: 
\begin{enumerate}[noitemsep]
    \item the formulation of a two-level hybrid restricted additive Schwarz preconditioner that uses the MS-GFEM approximation space based on the theory of \cite{Ma2025} as a coarse correction for con\-vec\-tion-diffusion problems;
    \item extension of the continuous MS-GFEM theory \cite{Ma2025} to convection-diffusion equations in conservative form, as well as a convergence proof of the iterative MS-GFEM for divergence-free velocity fields; 
    \item comprehensive numerical experiments in two and three dimensions for a variety of different discretization methods, demonstrating robustness with respect to P\'{e}clet number, diffusion heterogeneity, and the number of subdomains, including large-scale tests with more than $10^5$ subdomains; and
    \item demonstration of the method’s performance in the vanishing-diffusion limit (formally infinite P\'{e}clet number).
\end{enumerate}
The remainder of the paper is structured as follows: After stating the governing equations and introducing notation in \Cref{sec:model}, we formulate MS-GFEM for convection-diffusion problems in \Cref{sec:msgfem} and provide local and global error estimates on a continuous level. Subsequently, an iterative version of MS-GFEM is presented, which motivates the definition of the two-level hybrid restricted additive Schwarz preconditioner. This is detailed in \Cref{sec:discretization} in a generic discrete setting. In \Cref{sec:experiments}, we perform a thorough numerical investigation of the proposed preconditioner within the generalized minimal residual method (GMRES). 


\section{Governing equations}
\label{sec:model}

Let $\domain \subset \R^\dimDomain$, $\dimDomain \in \{2,3\}$, be a bounded domain with Lipschitz boundary 
$\boundary = \overline{\bdryD} \cup \overline{\bdryO}$, which is separated into a Dirichlet part, $\bdryD$, and 
an outflow part, $\bdryO$, with $\bdryD \cap \bdryO = \emptyset$.
We consider the stationary convection-diffusion equation in conservation form: Find $\fwp\colon\domain\rightarrow\R$ satisfying
\begin{subequations}
\label{eq:fwp}
\begin{alignat}{2}
  \label{eq:pde}
  \divergence(-\diffTensor \nabla \fwp + \velField \fwp) &= \rhssource &&\quad\quad\text{ in } \Omega,\\
  \label{eq:bcD}
  \fwp &= \bcD &&\quad\quad\text{ on } \bdryD,\\
  \label{eq:bcO}
  -\diffTensor \nabla \fwp \cdot \bn &= \bcO &&\quad\quad\text{ on } \bdryO,
\end{alignat}
\end{subequations}
where $\rhssource \in \Hone\dual$ is the source term, $\bn$ denotes the unit outer normal vector on $\boundary$, and $\bcD \in H^{1/2}(\bdryD)$, $\bcO \in H^{-1/2}(\bdryO)$ are the Dirichlet and outflow boundary conditions, respectively. Furthermore, $\velField \in \sobolevD{1}{\infty}$ is the velocity field satisfying $\bdryO \subseteq \{\velField \cdot \bn \geq 0\}$, and $\diffTensor \in \Lpdd{\infty}$ is the (pointwise) symmetric and uniformly elliptic diffusion tensor, i.e.,
\begin{equation}\label{eq:diffBounds}
 \diffMin \abs{\boldvar{\xi}}^2 \leq \diffTensor(\x) \boldvar{\xi} \cdot \boldvar{\xi} \leq \diffMax \abs{\boldvar{\xi}}^2 
 \quad \text{ for all } \x \in \Omega, \; \boldvar{\xi} \in \Rd,
\end{equation}
for some $\diffMin$, $\diffMax > 0$. Here, $\abs{\cdot}$ denotes the Euclidean norm on $\Rd$. In this work, we are primarily interested in convection-dominated problems, i.e., those with a high \emph{P\'{e}clet number}
\begin{equation}\label{eq:peclet}
  \peclet_L := \frac{L\norm{\velField}_{\Lpd{\infty}}}{\diffMin},
\end{equation}
where $L$ is a characteristic length scale, e.g., the mesh size $L=h$ in the discretized setting. We then refer to $\peclet_L = \peclet_h$ as the grid P\'{e}clet number.

Let $\fwpTestSpace := \{ \fwpTest \in \Hone \; | \; \fwpTest = 0 \text{ on } \bdryD \}$ and
\begin{equation*}
  \bilWF{\fwp}{\fwpTest}{} := \int_\domain \diffTensor \nabla \fwp \cdot \nabla \fwpTest \dx 
    - \int_\domain \fwp \velField \cdot \nabla \fwpTest \dx 
    + \int_{\bdryO} \fwp \fwpTest \velField \cdot \bn \ds.
\end{equation*}
For the analysis of the multiscale method, we eliminate the inhomogeneous Dirichlet boundary condition from the ansatz space in the weak formulation by selecting some function $\fwpDBC \in \Hone$ with $\fwpDBC = \bcD$ on $\bdryD$ and defining
\begin{equation*}
  \rhsWF{\fwpTest}{} := \int_\domain \rhssource \fwpTest \dx 
    - \int_{\bdryO} \bcO \fwpTest \ds - \bilWF{\fwpDBC}{\fwpTest}{},
\end{equation*}
which leads to the following weak form of the governing equations:
\begin{equation}\label{eq:weakFormGlobal}
  \text{Find $\fwpSol \in \fwpTestSpace$ such that } 
  \bilWF{\fwpSol}{\fwpTest}{} = \rhsWF{\fwpTest}{}
  \text{ for all $\fwpTest \in \fwpTestSpace$.}
\end{equation}
The weak solution of \cref{eq:fwp} is given by $\fwpSol + \fwpDBC$. 
For later use, we introduce the notation 
\begin{align}
  \notag 
  \bilWF{\fwp}{\fwpTest}{\ovlpdomain} 
    &:= \int_\ovlpdomain \diffTensor \nabla \fwp \cdot \nabla \fwpTest \dx 
    - \int_\ovlpdomain \fwp \velField \cdot \nabla \fwpTest \dx 
    + \int_{\bdryO \cap \partial \ovlpdomain} \fwp \fwpTest \velField \cdot \bn \ds,\\
    \notag 
    \rhsWF{\fwpTest}{\ovlpdomain} 
    &:= \int_\ovlpdomain \rhssource \fwpTest \dx 
    - \int_{\bdryO \cap \partial \ovlpdomain} \bcO \fwpTest \ds
    - \bilWF{\fwpDBC}{\fwpTest}{\ovlpdomain},\\
    \label{eq:bilIP}
    \bilIP{\fwp}{\fwpTest}{\ovlpdomain} 
    &:= \int_\ovlpdomain \diffTensor \nabla \fwp \cdot \nabla \fwpTest \dx, 
\end{align}
as well as the space
\begin{equation*}
  \fwpTestD{\ovlpdomain} := \{ \fwpTest \in H^1(\ovlpdomain) \; | \; \fwpTest = 0 \text{ on } \bdryD \cap \partial \ovlpdomain \}
\end{equation*}
for arbitrary subdomains $\ovlpdomain \subseteq \domain$.

\begin{remark}
    Throughout this paper, we assume that the global problem as well as all occurring subdomain problems are well-posed. Since our main goal is the numerical investigation of the proposed method, this allows us to focus on its development and properties. In practice, this might lead to some further assumptions on the coefficients $\diffTensor$ and $\velField$ or, in the discretized setting, on the mesh size $h$, which we will not specify explicitly. Related discussions can be found, e.g., in \cite[Sections 11.1, 11.5, and A.7]{ToselliWidlund2005}.
\end{remark}


\section{MS-GFEM for convection-diffusion problems}
\label{sec:msgfem}

Since the method we propose is based on a multiscale spectral generalized finite element method, we first formulate the multiscale method and prove convergence including error bounds on a continuous level. For this purpose, we follow the framework of \cite{Ma2025}. Note that \cite[Section 6.1]{Ma2025} deals with convection-diffusion equations in a non-conservative form, which coincides with the conservative form \cref{eq:fwp} if $\divergence \velField \equiv 0$. Hence, we effectively solve different equations if the velocity field is not divergence-free and we thus obtain different constants in the theory below.

Let $\{\ovlpdomain_\ddindex\}_{\ddindex=1}^{\nDoms}$ be overlapping subdomains of $\domain$ such that $\cup_{\ddindex=1}^{\nDoms} \ovlpdomain_\ddindex = \domain$. We introduce a partition of unity $\{\pu_\ddindex\}_{\ddindex=1}^\nDoms$ subordinate to the open cover $\{\ovlpdomain_\ddindex\}_{\ddindex=1}^{\nDoms}$ that satisfies 
\begin{equation}\label{eq:pu}
\supp(\pu_\ddindex) \subseteq \overline{\ovlpdomain_\ddindex}, \quad 
0 \leq \pu_\ddindex \leq 1, \quad 
\sum_{\ddindex=1}^{\nDoms} \pu_\ddindex = 1, \quad
\pu_\ddindex \in \sobolev{1}{\infty}, \quad 
\abs{\nabla \pu_\ddindex} \leq \frac{\puC}{\diam(\ovlpdomain_\ddindex)}
\end{equation}
for some $\puC > 0$ and all $\ddindex = 1,\ldots,\nDoms$.
The general idea is to construct approximation spaces on these overlapping subdomains that are tailored to the specific structure of the PDE. Locally, the solution is split into a \emph{particular function} $\fwpiPart$ satisfying the global boundary conditions and a function lying in the carefully designed $\locDim_\ddindex$-dimensional local approximation space $\locApproxSpace$, which will be defined in \cref{eq:weakFormLocal,eq:locApproxSpace}, respectively. The global particular function $\fwpPart$ and the global approximation space $\globalApproxSpace$ are then created by glueing the components together using a partition of unity,
\begin{equation*}
\fwpPart := \sum_{\ddindex=1}^{\nDoms} \pu_\ddindex \fwpiPart, \qquad
\globalApproxSpace := \Bigg\{\sum_{\ddindex=1}^{\nDoms} \pu_\ddindex \phi_\ddindex : \phi_\ddindex \in \locApproxSpace \Bigg\},
\end{equation*}
and the generalized finite element approximation of \cref{eq:weakFormGlobal} reads as follows:
\begin{equation}\label{eq:weakFormGFEM}
  \text{Find $\fwpGFEM \in \fwpPart + \globalApproxSpace$ such that } 
  \bilWF{\fwpGFEM}{\fwpTest}{} = \rhsWF{\fwpTest}{}
  \text{ for all $\fwpTest \in \globalApproxSpace$.}
\end{equation}
Since the global error is bounded by the local approximation errors as stated in \cref{th:msgfem:global} at the end of this section, the heart of the method lies in the construction of the local approximations.

For the definition of the local approximation spaces, we need slightly larger subdomains to obtain significantly better approximation properties. These \emph{oversampling domains} $\{\osdomain_\ddindex\}_{\ddindex=1}^{\nDoms}$ satisfy $\ovlpdomain_\ddindex \subset \osdomain_\ddindex \subset \domain$. On each oversampling domain, we define the local particular function $\fwpiPart$ as a local solution to the weak form of the PDE:
\begin{equation}\label{eq:weakFormLocal}
  \text{Find $\fwpiPart \in \fwpTestSpaceLoc$ such that } 
  \bilWF{\fwpiPart}{\fwpTest}{\osdomain_\ddindex} = \rhsWF{\fwpTest}{\osdomain_\ddindex}
  \text{ for all $\fwpTest \in \fwpTestSpaceLoc$,}
\end{equation}
where
\begin{equation*}
  \fwpTestSpaceLoc := \big\{ \fwp \in H^{1}(\osdomain_\ddindex) : \fwp = 0 \text{ on } (\osbdry \cap \bdryD) \cup (\osbdry \cap \domain) \big\} \subseteq \fwpTestD{\osdomain_\ddindex}.
\end{equation*}
Note that the boundary conditions on the interior subdomain boundary are arbitrary and mainly need to be chosen to formulate a well-posed subdomain problem. Depending on the specific application, other boundary conditions for the particular solution might be better suited.

A crucial observation for the construction of the local approximation spaces is that the local approximation error $(\fwpSol|_{\osdomain_\ddindex} - \fwpiPart) \in \fwpTestD{\osdomain_\ddindex}$ satisfies
\begin{equation}\label{eq:harmonicError}
  \bilWF{\fwpSol|_{\osdomain_\ddindex} - \fwpiPart}{\fwpTest}{\osdomain_\ddindex} = 0 \quad
  \text{ for all $\fwpTest \in \fwpTestSpaceLoc$,}
\end{equation}
i.e., it lies in the problem-specific operator-harmonic (or $B$-harmonic) space
\begin{equation*}
  \harmSpace := \big\{ \fwp \in \fwpTestD{\osdomain_\ddindex} :  \bilWF{\fwp}{\fwpTest}{\osdomain_\ddindex} = 0 \quad \forall \fwpTest \in \fwpTestSpaceLoc\big\},
\end{equation*}
motivating to correct the error with functions from this space. As we will see later, the $B$-harmonic space can be approximated well by low-dimensional subspaces, which is key to the efficiency of our method. Note that $\harmSpace$ equipped with $\bilIP{\cdot}{\cdot}{\osdomain_\ddindex}$ from \cref{eq:bilIP} is a Hilbert space unless $\osbdry \cap \bdryD = \emptyset$ and $\divergence \velField \equiv 0$ on $\osdomain_\ddindex$, in which case positive-definiteness is lost due to the existence of (non-zero) constant functions in $\harmSpace$. In either case, the optimal local approximation spaces can be constructed by means of the same generalized eigenvalue problem as outlined below (see \cite[Section 2.2]{Ma2025} for the proof).

We start by defining the restriction operator
\begin{equation}\label{eq:resOp}
  \resOp_\ddindex \colon \harmSpace \to \fwpTestSpaceLoc[\ovlpdomain], \quad \fwp \mapsto \pu_\ddindex \fwp|_{\ovlpdomain_\ddindex},
\end{equation}
and we equip both the domain and codomain with the (semi-)norm induced by the elliptic part of the bilinear form, $\norm{\cdot}_{a,\ovlpdomain} := (\bilIP{\cdot}{\cdot}{\ovlpdomain})^{1/2}$. Approximating the image of $\resOp_\ddindex$ boils down to calculating a singular value decomposition, which leads to the following weak form of the local generalized eigenvalue problem: Find $\eval \in [0,+\infty]$\footnote{The case $\eval = +\infty$ corresponds to the existence of $\evec \in \harmSpace$ such that $\bilIP{\evec}{\psi}{\osdomain_\ddindex} = 0$ for all $\psi \in \harmSpace$ but $\bilIP{\pu_\ddindex\evec}{\pu_\ddindex\widehat{\psi}}{\ovlpdomain_\ddindex} \neq 0$ for some $\widehat{\psi} \in \harmSpace$.}, $\evec \in \harmSpace$ such that
\begin{equation}\label{eq:gevp}
  \bilIP{\pu_\ddindex\evec}{\pu_\ddindex\psi}{\osdomain_\ddindex}  = \eval \bilIP{\evec}{\psi}{\osdomain_\ddindex} \quad\text{for all $\psi \in \harmSpace$.}
\end{equation}
Note that, despite being based on the non-symmetric convection-diffusion equations, the above eigenvalue problem consists of two symmetric positive (semi-)definite bilinear forms, yielding only real and non-negative eigenvalues. The asymmetry of the underlying problem and in particular the velocity field itself are only incorporated through the harmonic constraint in the ansatz space of the eigenfunctions, $\harmSpace$. Furthermore, the kernels of the bilinear forms on both sides of \cref{eq:gevp} have a trivial intersection due to the partition of unity and the discretized eigenproblem is thus non-defective. From a computational point of view, iterative solvers for this eigenproblem typically converge fast since the relative spectral gap between consecutive eigenvalues is large (cf.\ \cref{th:msgfem:evaldecay} below, where an exponential decay is shown).

Denoting the $\evalindex$-th eigenpair of \cref{eq:gevp} as $(\evalloc, \evecloc)$, where the eigenvalues $\evalloc[1] \geq \evalloc[2] \geq \ldots$ are numbered in non-ascending order, the local approximation space $\locApproxSpace$ is built from the eigenfunctions corresponding to the $\locDim_\ddindex$ largest eigenvalues:
\begin{equation}\label{eq:locApproxSpace}
  \locApproxSpace := \mathrm{span} \big\{ \evecloc[1]|_{\ovlpdomain_\ddindex}, \ldots, \evecloc[\locDim_\ddindex]|_{\ovlpdomain_\ddindex} \big\}.
\end{equation}

\begin{remark}\label{rm:gevp}
Solving the (discretized) local eigenproblem \cref{eq:gevp} requires calculating a basis for the local $B$-harmonic space. This involves solving as many local boundary value problems as there are degrees of freedom in the discretization of $\partial \osdomain_\ddindex$, which is particularly expensive on finer meshes. Hence, we avoid the explicit calculation of a basis for $\harmSpace$ by introducing a Lagrange multiplier and enforcing the $B$-harmonic constraint weakly: Find $\eval \in [0,+\infty]$, $\evec \in \fwpTestD{\osdomain_\ddindex}$, $\zeta \in \fwpTestSpaceLoc$ such that
\begin{subequations}\label{eq:gevpSaddlePoint}
\begin{alignat}{4}
  \eval\big(\bilIP{\evec}{\psi}{\osdomain_\ddindex} &+ \emph{$\bilWF{\psi}{\zeta}{\osdomain_\ddindex}$}\big) &&= \bilIP{\pu_\ddindex\evec}{\pu_\ddindex\psi}{\osdomain_\ddindex} \quad&&\text{for all $\psi \in \fwpTestD{\osdomain_\ddindex}$},&&\\
  &\phantom{{}+{}}\bilWF{\evec}{\xi}{\osdomain_\ddindex} &&= 0\quad&&\text{for all $\xi \in \fwpTestSpaceLoc$.}&&
\end{alignat}
\end{subequations}
Although this alternative formulation of the eigenproblem has almost twice as many unknowns as the original problem, it is not only more convenient for the implementation due to its formulation on standard spaces, but also computationally cheaper. Note that the eigenfunctions of \cref{eq:gevpSaddlePoint} corresponding to non-zero eigenvalues coincide with the eigenfunctions of \cref{eq:gevp}. Therefore, both formulations are equivalent for creating the local approximation spaces and we employ the mixed formulation \cref{eq:gevpSaddlePoint} in our numerical simulations.
\end{remark}

Two central assumptions for the MS-GFEM theory are a Caccioppoli-type inequality as well as a weak approximation property, which we prove next.

\begin{lemma}[Caccioppoli-type inequality]\label{th:caccioppoli}
    Let $D \subset D^* \subset \domain$ and let $\eta \in C^1(\overline{D^*})$ be such that $\eta = 0$ on $\partial D^* \cap \domain$. 
    Then, the following holds for all $\fwp$, $\fwpTest \in H_B(D^*)$:
    \begin{equation}\label{eq:cacId}
        \bilIP{\eta\fwp}{\eta\fwpTest}{D^*} 
        = \int_{D^*} \left(\diffTensor \nabla \eta \cdot \nabla \eta\right)\fwp \fwpTest\dx 
        - \int_{D^*} \eta^2 \fwp\fwpTest \divergence \velField \dx
        - \frac{1}{2}\int_{D^*} \eta^2 \velField \cdot \nabla(\fwp\fwpTest)\dx.
    \end{equation}
    If in addition $\delta := \mathrm{dist}(D, \partial D^*\setminus\boundary) >0$, then    
    \begin{equation}\label{eq:cac}
        \norm{\fwp}_{a,D} \leq 
        C_{\mathrm{cac}}^{\mathrm{I}} \delta^{-1} \norm{\fwp}_{L^2(D^*\setminus D)} 
        + C_{\mathrm{cac}}^{\mathrm{II}}\norm{\fwp}_{L^2(D^*)} \quad \text{for all} \quad \fwp \in H_B(D^*),
    \end{equation}
    where
    \begin{align*}
        C_{\mathrm{cac}}^{\mathrm{I}} &:= \sqrt{2\diffMax + 1} \quad \text{and} \quad
        C_{\mathrm{cac}}^{\mathrm{II}} &:= \sqrt{2\norm{\divergence \velField}_{\Lp{\infty}} + \norm{\velField}^2_{\Lpd{\infty}}\left(1 + \diffMin^{-1} \right)}.
    \end{align*}
\end{lemma}

\begin{proof}
    Let $\fwp$, $\fwpTest \in H_B(D^*)$ be arbitrary. Applying \cite[Proposition 5.3]{Ma2025} with $b_{\textbf{L}}(\vec{w}_1, \vec{w}_2) := \int_{D^*} \diffTensor\, \vec{w}_1 \cdot \vec{w}_2 \dx$ and the definitions of \cite[Example 5.1]{Ma2025} yields
    \begin{equation}\label{eq:proof:cac1}
        \begin{aligned}
            \bilIP{\eta\fwp}{\eta\fwpTest}{D^*} 
            &= \int_{D^*} \left(\diffTensor \nabla \eta \cdot \nabla \eta\right)\fwp \fwpTest\dx \\
            &\phantom{={}} + \frac{1}{2} \left( 
            \int_{D^*} \diffTensor \nabla \fwp \cdot \nabla (\eta^2 \fwpTest)\dx 
            + \int_{D^*} \diffTensor \nabla \fwpTest \cdot \nabla (\eta^2 \fwp)\dx 
            \right).
        \end{aligned}
    \end{equation}
    Since $\eta^2\fwpTest \in H^1_{D,0}(D^*)$, the harmonic constraint on $\fwp$ yields
    \begin{equation*}
        \begin{aligned}
            \int_{D^*} \diffTensor \nabla \fwp \cdot \nabla (\eta^2 \fwpTest)\dx 
            &= \int_{D^*} \fwp \velField \cdot \nabla (\eta^2 \fwpTest)\dx 
            - \int_{\bdryO \cap \partial D^*} \eta^2 \fwp \fwpTest \velField \cdot \bn \ds\\
            &= - \int_{D^*} \eta^2 \fwpTest \divergence (\fwp \velField) \dx\\
            &= - \int_{D^*} \eta^2 \fwp \fwpTest \divergence \velField \dx
            - \int_{D^*} \eta^2 \fwpTest \velField \cdot \nabla \fwp \dx,
        \end{aligned}
    \end{equation*}
    where we integrated by parts and used the boundary conditions in the second line as well as the product rule in the third. Switching the roles of $\fwp$ and $\fwpTest$, combining the resulting equations with \cref{eq:proof:cac1}, and applying the product rule then yields \cref{eq:cacId}.

    To prove the Caccioppoli-type inequality \cref{eq:cac}, let $\fwp = \fwpTest$ in \cref{eq:cacId} to obtain
    \begin{equation}\label{eq:proof:cac2}
        \norm{\eta\fwp}^2_{a,D^*} 
        = \int_{D^*} \left(\diffTensor \nabla \eta \cdot \nabla \eta\right)\fwp^2 \dx 
        - \int_{D^*} \eta^2 \fwp^2 \divergence \velField \dx
        + \int_{D^*} \eta \fwp \velField \cdot (u \nabla \eta - \nabla(\eta\fwp))\dx.
    \end{equation}
    %
    Using (pointwise) weighted Young's inequalities, the integrand of the third integral in the above equation can be bounded from above by 
    \begin{equation*}
        \frac{\abs{\eta\fwp\velField}^2}{2} + \frac{\abs{\fwp\nabla\eta}^2}{2} + \frac{\abs{\eta\fwp\velField}^2}{2\diffMin} + \frac{\diffMin\abs{\nabla(\eta\fwp)}^2}{2},
    \end{equation*}
    such that H\"{o}lder's inequality and the spectral bound \cref{eq:diffBounds} yield the following upper bound for this integral:
    \begin{equation*}
    \begin{aligned}
        \frac{1}{2}\Big(&\norm{\eta\fwp}^2_{L^2(D^*)}\norm{\velField}^2_{\Lpd{\infty}}
        + \norm{\fwp\nabla\eta}^2_{L^2(D^*;\,\R^\dimDomain)}\\
        &+ \diffMin^{-1}\norm{\eta\fwp}^2_{L^2(D^*)}\norm{\velField}^2_{\Lpd{\infty}}
        + \norm{\eta\fwp}^2_{a,D^*}
        \Big).
    \end{aligned}
    \end{equation*}
    Applying H\"{o}lder's inequality to the first and second integral of \cref{eq:proof:cac2} and combining the estimates yields
    \begin{equation*}
        \begin{aligned}
            \norm{\eta\fwp}^2_{a,D^*} &\leq
            \norm{\fwp\nabla\eta}^2_{L^2(D^*;\,\R^\dimDomain)} \left(2\diffMax + 1\right)\\
            &\phantom{\leq{}}+ \norm{\eta\fwp}^2_{L^2(D^*)} \left[2\norm{\divergence\velField}_{\Lp{\infty}} + \norm{\velField}^2_{\Lpd{\infty}}(1 + \diffMin^{-1}) \right].
        \end{aligned}
    \end{equation*}
    Choosing a cut-off function $\eta \in C^1(\overline{D^*})$ such that $0\leq\eta\leq1$ and $\abs{\nabla\eta} \leq \delta^{-1}$ in $D^*$, and $\eta = 1$ in $D$, we obtain
    \begin{equation*}
    \begin{aligned}
        \norm{\fwp}^2_{a,D} 
        \leq \norm{\eta\fwp}^2_{a,D^*}
        &\leq (2\diffMax + 1) \delta^{-2} \norm{\fwp}^2_{L^2(D^*\setminus D)}\\
        &\phantom{\leq{}}+ \left(2\norm{\divergence\velField}_{\Lp{\infty}} + \norm{\velField}^2_{\Lpd{\infty}}(1 + \diffMin^{-1})\right)
        \norm{\fwp}^2_{L^2(D^*)}.
    \end{aligned}
    \end{equation*}
    Taking the square root establishes \cref{eq:cac}.
\end{proof}

\begin{lemma}[Weak approximation property]\label{th:weakApprox}
    Let $D \subset D^* \subset D^{**} \subset \Omega$ 
    such that $\delta := \mathrm{dist}(D^*, \partial D^{**}\setminus\boundary) >0$ and $\mathrm{dist}(\x, \partial D^*\setminus\boundary) \leq \delta$ for all $\x \in \partial D \setminus \boundary$, and let $V_\delta(D^* \setminus D) := \big\{\x \in D^{**} : \mathrm{dist}(\x, D^* \setminus D) \leq \delta \big\}$ be the $\delta$-vicinity of the set $D^* \setminus D$. 
    
    Then, there are constants $C_1$, $C_2 > 0$ depending only on $\dimDomain$ such that for every integer $N \geq C_1 \abs{D^{**}}\delta^{-\dimDomain}$ and every $a$, $b > 0$, there exists an $N$-dimensional space $\mathcal{Q}_N(D^{**}) \subset L^2(D^{**})$ such that the following estimate holds for all $\fwp \in H_B(D^{**})$:
    \begin{equation*}
    \begin{aligned}
        \inf_{\fwpTest \in \mathcal{Q}_N(D^{**})} &\left(a \norm{\fwp - \fwpTest}_{L^2(D^*\setminus D)} + b \norm{\fwp - \fwpTest}_{L^2(D^*)} \right)\\
        &\qquad \leq C_2 \, \diffMin^{-1/2} \left(a \abs{V_\delta(D^* \setminus D)}^{1/\dimDomain} + b \abs{D^{**}}^{1/\dimDomain} \right) N^{-1/\dimDomain} \norm{\fwp}_{a,D^{**}}.
    \end{aligned}
    \end{equation*}
\end{lemma}

\begin{proof}
   Due to \cref{eq:diffBounds}, we have $\norm{\nabla \fwp}_{L^2(D^{**})} \leq \diffMin^{-1/2} \norm{\fwp}_{a,D^{**}}$. Thus, the statement follows from case (iii) of \cite[Lemma 5.5]{Ma2025}.
\end{proof}

With the above preparations, we obtain the following local error estimate:
\begin{theorem}[Local error estimate]\label{th:msgfem:local}
Let $\fwpSol$, $\fwpiPart$, and $\locApproxSpace$ be defined via \cref{eq:weakFormGlobal,eq:weakFormLocal,eq:locApproxSpace}, respectively. Then, for each $\ddindex=1,\ldots,\nDoms$, we have
\begin{equation}\label{eq:msgfem:local}
  \inf_{\evec \in \fwpiPart|_{\ovlpdomain_\ddindex} + \locApproxSpace} 
  \norm{\pu_\ddindex(\fwpSol|_{\ovlpdomain_\ddindex} - \evec)}_{a,\ovlpdomain_\ddindex} 
  \leq \sqrt{\evalloc[\locDim_{\ddindex+1}]}
  \norm{\fwpSol|_{\osdomain_\ddindex} - \fwpiPart}_{a,\osdomain_\ddindex}.
\end{equation}
\end{theorem}
\begin{proof}
The statement follows from \cite[Theorem 2.23]{Ma2025}: Assumption 2.13 of \cite{Ma2025} is satisfied with the bilinear form $\bilIP{\cdot}{\cdot}{D}$ from \cref{eq:bilIP}, for every domain $\ovlpdomain_\ddindex \subset D \subset \osdomain_\ddindex$. In particular, if $\partial D$ does not intersect the global Dirichlet boundary and $\divergence \velField \equiv 0$ on $D$, then the space $\mathcal{K}_D$ in \cite[Assumption 2.13]{Ma2025} is the one-dimensional space of constant functions. Otherwise, $\bilIP{\cdot}{\cdot}{D}$ defines an inner product on $H_B(D)$. Assumption 2.17 of \cite{Ma2025} follows from \cite[Proposition 3.7]{Ma2025} due to the Caccioppoli-type inequality (\cref{th:caccioppoli}) and the weak approximation property (\cref{th:weakApprox}). Note that the Kolmogorov $n$-width given in \cite[Theorem 2.23]{Ma2025} coincides with $\sqrt{\evalloc[\locDim_{\ddindex+1}]}$ due to \cite[Lemmas 2.19 and 2.20]{Ma2025}.
\end{proof}

Crucial for the efficiency of the method is a low-dimensional approximation property of $\harmSpace$ by $\locApproxSpace$. To keep the presentation simple, we assume that the overlapping domains $\ovlpdomain_\ddindex$ and oversampling domains $\osdomain_\ddindex$ are truncated concentric cubes with side lengths $H_\ddindex$ and $H_\ddindex^*$, respectively, and we define $\delta_\ddindex^* := H_\ddindex^* - H_\ddindex > 0$. For a discussion of more general domains, we refer to \cite[Remarks 3.11 and 3.15]{Ma2025}.
\begin{theorem}[Exponential convergence of local approximations]\label{th:msgfem:evaldecay}
There are constants $C_\ddindex$, $c_\ddindex > 0$ and $N_\ddindex \in \N$ such that
\begin{equation}\label{eq:msgfem:evaldecay}
  \sqrt{\evalloc[N]} \leq C_\ddindex e^{-c_\ddindex N^{1/\dimDomain}}
  \quad \text{for all $N > N_\ddindex$}.
  \end{equation}
\end{theorem}
\begin{proof}
    The nearly exponential decay of the eigenvalues as well as the constants given in \cref{rm:constants} below are consequences of \cite[Theorem 3.8 (iii)]{Ma2025}. The validity of 
    \cite[Assumption 2.13]{Ma2025} has already been discussed above, while \cite[Assumptions 3.1 and 3.4]{Ma2025} are satisfied due to \cref{th:caccioppoli,th:weakApprox}, respectively. 
\end{proof}

\begin{remark}\label{rm:constants}
Defining 
\begin{align*}
  \Theta_\ddindex &= \widehat{C} \left(\frac{2\diffMax + 1}{\diffMin} \right)^{1/2} 
  \bigg(\frac{H_\ddindex^*}{\delta_\ddindex^*} \bigg)^{(\dimDomain -1)/\dimDomain},\\ 
  \sigma_\ddindex^* &= \frac{\delta_\ddindex^*}{4} \left(\frac{2\norm{\divergence \velField}_{\Lp{\infty}} + \norm{\velField}^2_{\Lpd{\infty}}\left(1 + \diffMin^{-1} \right)}{2\diffMax + 1} \right)^{1/2}, 
\end{align*}
the constants in \cref{th:msgfem:evaldecay} can be quantified:
\begin{align*}
  c_\ddindex &= \left(2e\Theta_\ddindex + 1\right)^{-1},\\
  N_\ddindex &= 2\max\Bigg\{(3\dimDomain)^{-1}(2e\Theta_\ddindex)^\dimDomain{\sigma_\ddindex^*}^{\tfrac{\dimDomain^2}{\dimDomain-1}}\bigg(\frac{H_\ddindex^*}{\delta_\ddindex^*}\bigg)^{\dimDomain/(\dimDomain-1)},
  \left(2e\Theta_\ddindex\right)^{-\tfrac{\dimDomain}{\dimDomain-1}} \Bigg\},\\
  C_\ddindex &= \widetilde{C} \Bigg[\peclet_{H_\ddindex}
  \bigg(1 + \frac{\diffMin \norm{\divergence \velField}_{\Lp{\infty}}}{\norm{\velField}^2_{\Lpd{\infty}}}\bigg)^{1/2} 
  + \left(\frac{\diffMax}{\diffMin}\right)^{1/2}
  \Bigg],
\end{align*}
where $\widehat{C}$, $\widetilde{C} > 0$ are constants that are independent of problem-specific parameters. This provides us with two switches to control the local error \cref{eq:msgfem:local}: First, we can increase the number of eigenfunctions in the local approximation spaces. Second, increasing the oversampling size also contributes to a faster (and earlier) decay of the eigenvalues. Both factors can be observed in the numerical experiments in \Cref{sec:experiments}.
Another interesting aspect is that the slope $c_\ddindex$ of the nearly exponential decay \cref{eq:msgfem:evaldecay} does not depend on the velocity field, whereas $C_\ddindex$ and $N_\ddindex$  do. More precisely, if the velocity field is divergence-free, the first term of $C_\ddindex$ reduces to the P\'{e}clet number $\peclet_{H_\ddindex}$ associated to the coarse mesh. While this suggests that larger approximation spaces $\locApproxSpace$ are needed the more convection-dominated the problems become, we only observe a mild influence in our numerical results; the method is indeed very robust with respect to changes in problem parameters when used as a preconditioner.
\end{remark}

\begin{remark}\label{rm:locBC}
As pointed out before, the local approximation error \cref{eq:msgfem:local} can potentially be improved by choosing more appropriate local boundary conditions for the particular functions \cref{eq:weakFormLocal}. However, due to the nearly exponential decay of the eigenvalues and our use of the method as a preconditioner, the influence of adjusting the local approximation spaces by selecting (only a few) more eigenfunctions is expected to be much bigger.
\end{remark}

For the global error estimate we need some additional notation. Let
\begin{equation}\label{eq:colConst}
  \colConst := \min \Big\{k \in \N \, : \, \mathrm{card}\{\ddindex : \x \in \ovlpdomain_\ddindex\} \leq k \; \forall\; \x \in \domain \Big\}
\end{equation}
be the coloring constant of the overlapping subdomains $\{\ovlpdomain_\ddindex\}_{\ddindex=1}^{\nDoms}$, 
\begin{equation}\label{eq:evalMax}
  \evalMax := \max_{\ddindex=1,\ldots,\nDoms} \evalloc[\locDim_{\ddindex+1}]
\end{equation}
the largest eigenvalue whose eigenfunction is not included in the coarse space and
\begin{equation*}
  H^*_{\mathrm{max}} := \max_{\ddindex=1,\ldots,\nDoms} \mathrm{diam}(\osdomain_\ddindex)
\end{equation*}
the maximal diameter of the oversampling subdomains. If the global bilinear form $\bilB$ is coercive (e.g., when $\divergence \velField \equiv 0$), \cref{th:msgfem:global} below always holds. In the general indefinite case, resolution conditions on $\evalMax$ and $H^*_{\mathrm{max}}$ are required. While these conditions can be explicitly stated, similar to the constants in \cref{th:msgfem:evaldecay}, the resulting estimates from \cite[Corollary 3.27]{Ma2025} are overly pessimistic. See \cite[Sections 6.1--6.2]{Ma2025} for a precise statement of these conditions in a special case, as well as discussions of related problems. Since we do not observe this influence on the preconditioner in the numerical experiments of \Cref{sec:experiments}, we only state a simplified form of the global error estimate:

\begin{theorem}[Quasi-optimal global convergence]\label{th:msgfem:global}
Let $\fwpSol$ be the solution of \cref{eq:weakFormGlobal} and $\fwpGFEM$ the MS-GFEM approximation defined by \cref{eq:weakFormGFEM}. 
If  $\bilWF{\cdot}{\cdot}{}$ is coercive or $\evalMax$ and $H^*_{\mathrm{max}}$ are sufficiently small, then
\begin{equation}\label{eq:msgfem:global}
  \norm{\fwpSol - \fwpGFEM}_{a,\domain} \leq C_B \bigg(\colConst \, \evalMax \sum_{\ddindex=1}^M \norm{\fwpSol|_{\osdomain_\ddindex} - \fwpiPart}_{a,\osdomain_\ddindex}^2 \bigg)^{1/2},
\end{equation}
where $C_B \sim 1 + \diffMin^{-1}\left(\norm{\divergence \velField}_{\Lp{\infty}} + \norm{\velField}_{\Lpd{\infty}}\right)$ is the operator norm of the bilinear form $\bilWF{\cdot}{\cdot}{}$ up to coefficient-independent constants.
\end{theorem}
\begin{proof}
    For coercive $\bilB$ the estimate follows immediately from \cite[Theorem 3.20]{Ma2025}. In the general indefinite case, the theorem is a consequence of \cite[Corollary 3.27]{Ma2025}, where \cite[Assumptions 3.23 and 3.24]{Ma2025} are proved like \cite[Lemmas 6.4 and 6.5]{Ma2025}. 
\end{proof}

\begin{corollary}
If $\divergence \velField \equiv 0$ in $\domain$, the estimate in Theorem \ref{th:msgfem:global} simplifies to
        \begin{equation}\label{eq:msgfem:global:coercive}
          \norm{\fwpSol - \fwpGFEM}_{a,\domain} \leq 
          \left(1 + C\,\peclet_H\right)^2 
          \sqrt{\colConst\colConst^* \evalMax }
          \norm{\fwpSol}_{a,\domain},
        \end{equation}
        where $C := \max_{\ddindex=1,\ldots,\nDoms} C_{\osdomain_\ddindex}$ with the Poincar\'{e} constants $C_{\osdomain_\ddindex}$ of $\osdomain_\ddindex$  as in \cite[Corollary A.15]{ToselliWidlund2005}, $\peclet_H := \max_{\ddindex=1,\ldots,\nDoms} \peclet_{H_\ddindex^*}$, and $\colConst^*$ is the coloring constant of the oversampling subdomains $\{\osdomain_\ddindex\}_{\ddindex=1}^{\nDoms}$ in analogy to \cref{eq:colConst}.
\end{corollary}
\begin{proof}
    First, observe that due to \cref{eq:harmonicError}, since $\fwpiPart \in \fwpTestSpaceLoc$, we have
    \begin{equation*}
        \bilWF{\fwpSol - \fwpiPart}{\fwpiPart}{\osdomain_\ddindex} = 0.
    \end{equation*}
    Using the coercivity and continuity of the bilinear form, we obtain
    \begin{align*}
        \norm{\fwpSol - \fwpiPart}_{a,\osdomain_\ddindex}^2 &\leq \abs{\bilWF{\fwpSol - \fwpiPart}{\fwpSol - \fwpiPart}{\osdomain_\ddindex}}
        = \abs{\bilWF{\fwpSol - \fwpiPart}{\fwpSol}{\osdomain_\ddindex}}\\
        &\leq
        \left(1 + C_{\osdomain_\ddindex} H_\ddindex^* \diffMin^{-1} \norm{\velField}_{\Lpd{\infty}}\right)
        \norm{\fwpSol - \fwpiPart}_{a,\osdomain_\ddindex}
        \norm{\fwpSol}_{a,\osdomain_\ddindex}.
    \end{align*}
    %
   Estimate \cref{eq:msgfem:global:coercive} then follows from \cref{eq:msgfem:global} and the definition of $\colConst^*$.
\end{proof}

\subsection{Iterative MS-GFEM}
\label{sec:iterative}

In \Cref{sec:experiments}, we employ MS-GFEM as a preconditioner within GMRES. For this purpose, we first formulate an iterative version of MS-GFEM that motivates the choice of preconditioner. This perspective has previously been taken in the context of elliptic \cite{Strehlow2024} and Helmholtz problems \cite{Ma2025Helmholtz}. 

The multiscale method presented in the last section can naturally be seen as a domain decomposition method: The local particular functions $\fwpiPart$ take the role of subdomain solves on the oversampling domains $\osdomain_\ddindex$, which are then prolongated via multiplication by the partition of unity $\pu_\ddindex$. On the other hand, the global approximation space $\globalApproxSpace$ can be seen as a coarse space in the following sense: Denoting $\fwpGFEM = \fwpPart + \fwp^s$, it follows from \cref{eq:weakFormGFEM} that $\fwp^s \in \globalApproxSpace$ is the unique solution of 
\begin{equation*}
  \bilWF{\fwp^s}{\fwpTest}{} = \bilWF{\fwpSol - \fwpPart}{\fwpTest}{} \quad \forall \; \fwpTest \in \globalApproxSpace,
\end{equation*}
i.e., the MS-GFEM solution $\fwpGFEM$ is a coarse space correction of the one-level restricted additive Schwarz approximation $\fwpPart$. Repeating this procedure of local subdomain solves and subsequent coarse space corrections using local Dirichlet data from the previous iteration defines an iterative version of MS-GFEM.

To state the method precisely, we first define the local projection $\projB\colon \fwpTestSpace \to \fwpTestSpaceLoc$ via
\begin{equation*}
  \bilWF{\projB(\fwpTest)}{w}{\osdomain_\ddindex} = \bilWF{\fwpTest|_{\osdomain_\ddindex}}{w}{\osdomain_\ddindex} \quad \forall \; w \in \fwpTestSpaceLoc,
\end{equation*}
and the coarse space projection $\projCoarse\colon \fwpTestSpace \to \globalApproxSpace$ via
\begin{equation*}
  \bilWF{\projCoarse(\fwpTest)}{w}{} = \bilWF{\fwpTest}{w}{} \quad \forall \; w \in \globalApproxSpace.
\end{equation*}
Additionally, we denote by $\msgfemOperator$ the linear operator, which maps a function $w \in \fwpTestSpace$ to the MS-GFEM solution $\fwpGFEM$ from \cref{eq:weakFormGFEM} with right-hand side $\rhsF := \bilWF{w}{\cdot}{}$. Following the explanation of the previous paragraph, this operator can be written as
\begin{equation*}
  \msgfemOperator w = \sum_{\ddindex = 1}^\nDoms \pu_\ddindex \projB(w)
  + \projCoarse\left(w - \sum_{\ddindex = 1}^\nDoms \pu_\ddindex \projB(w)\right).
\end{equation*}
Given some initial guess $\fwp^{(0)} \in \fwpTestSpace$, the MS-GFEM iteration is then defined as
\begin{equation}\label{eq:msgfem:iteration}
  \fwp^{(k+1)} := \fwp^{(k)} + \msgfemOperator(\fwpSol - \fwp^{(k)}), \quad k=0,1,\ldots,
\end{equation}
which can be seen as a preconditioned Richardson iteration. 
Note that the exact solution $\fwpSol$, which occurs in \cref{eq:msgfem:iteration}, is not needed to execute the iteration due to the relation $\bilWF{\fwpSol}{\fwpTest}{} = \rhsWF{\fwpTest}{}$; see \cref{eq:msgfem:iteration:discrete} for the corresponding discrete iteration.

In the following, we state a convergence result of the MS-GFEM iteration for the important special case of divergence-free velocity fields, which occurs, e.g., for incompressible fluids. Note that the proof remains the same for general coercive $\bilB$; we simply focus on divergence-free velocity fields for ease of presentation. 

\begin{theorem}[Convergence of the MS-GFEM iteration]\label{th:msgfem:iteration}
    Let $\divergence \velField \equiv 0$ in $\domain$ and let $\fwpSol$ be the solution of \cref{eq:weakFormGlobal}. 
    Then, for any $\vartheta \in (0,1)$, it is possible to choose the dimensions $\locDim_1, \ldots, \locDim_\nDoms \in \N$ of the local approximation spaces $\locApproxSpace$ such that
        \begin{equation}\label{eq:msgfem:iteration:1}
          \norm{\fwpTest - \msgfemOperator\fwpTest}_{a,\domain} \leq \vartheta
          \norm{\fwpTest}_{a,\domain} \quad \text{for all $\fwpTest \in \fwpTestSpace$,}
        \end{equation}
       and the sequence of MS-GFEM iterates $\{\fwp^{(k)}\}_{k=0}^\infty$ defined by \cref{eq:msgfem:iteration} converges for any initial guess $\fwp^{(0)} \in \fwpTestSpace$ strictly monotonically:
        \begin{equation}\label{eq:msgfem:iteration:2}
          \norm{\fwp^{(k+1)} - \fwpSol}_{a,\domain} \leq \vartheta
          \norm{\fwp^{(k)} - \fwpSol}_{a,\domain}, \quad k=0,1,\ldots
        \end{equation}
\end{theorem}

\begin{proof}
    Let $\vartheta \in (0,1)$ and $\fwpTest \in \fwpTestSpace$ be arbitrary. As mentioned above, $\msgfemOperator\fwpTest$ is the MS-GFEM approximation of the weak form \cref{eq:weakFormGlobal} with right-hand side $\rhsF := \bilWF{\fwpTest}{\cdot}{}$, which has the exact solution $\fwpTest$. Hence, the existence of $\locDim_1, \ldots, \locDim_\nDoms \in \N$ such that \cref{eq:msgfem:iteration:1} holds follows from \cref{eq:evalMax}, \cref{eq:msgfem:global:coercive} and \cref{th:msgfem:evaldecay}. Subtracting $\fwpSol$ from both sides of \cref{eq:msgfem:iteration}, we see that
    \begin{equation*}
        \fwp^{(k+1)} - \fwpSol = (\fwp^{(k)} - \fwpSol) - \msgfemOperator(\fwp^{(k)} - \fwpSol),
    \end{equation*}
    and estimate \cref{eq:msgfem:iteration:1} thus implies \cref{eq:msgfem:iteration:2}.
\end{proof}


\section{Discretization}
\label{sec:discretization}

So far, we have formulated MS-GFEM and the corresponding iteration on a continuous level. Since the method only requires the ability to perform subdomain solves and the solution of generalized eigenvalue problems on the same subdomains, it can easily be applied for different discretizations. Hence, we formulate the preconditioner we employ in \Cref{sec:experiments} in a generic discrete setting. 

The mesh of the computational domain is denoted by $\mesh$ and assumed to be sufficiently fine to resolve the fine-scale features of the coefficients. We further assume that the diffusion tensor is constant on each cell of the mesh. In the discrete setting, we only consider cases in which the overlapping subdomains $\ovlpdomain_\ddindex$ and the oversampling subdomains $\osdomain_\ddindex$ are a union of mesh elements of $\mesh$.

\subsection{Two-level hybrid restricted additive Schwarz preconditioner}
\label{sec:ras}

Let $\discRestriction\transpose$ be the matrix representation of the discretized extension-by-zero operator from $\fwpTestSpaceLoc$ to $\fwpTestSpace$, denote the matrix representation of the discretized embedding of $\globalApproxSpace$ into $\fwpTestSpace$ by $\discCoarseRestriction\transpose$, and let $\discBilForm$ be the stiffness matrix corresponding to the bilinear form $\bilB$. Defining $\discBilForm_\ddindex := \discRestriction\discBilForm\discRestriction\transpose$ and $\discBilForm_S := \discCoarseRestriction\discBilForm\discCoarseRestriction\transpose$, the matrix forms of the projections $\projB$ and $\projCoarse$ are
\begin{equation*}
  \discProjB = \discBilForm_\ddindex^{-1}\discRestriction\discBilForm
  \quad \text{and} \quad 
  \discProjCoarse = \discBilForm_S^{-1}\discCoarseRestriction\discBilForm,
\end{equation*}
respectively. The matrix representation of the discretized MS-GFEM map $\msgfemOperator$ can thus be written as $\mat{\msgfemOperator} = \discPreconditioner\discBilForm$ with the preconditioner
\begin{equation*}
  \discPreconditioner := 
  \sum_{\ddindex=1}^\nDoms \discRestriction\transpose \discPU \discBilForm_\ddindex^{-1}\discRestriction
  + \discCoarseRestriction\transpose\discBilForm_S^{-1}\discCoarseRestriction
  \left(\mat{I} - \sum_{\ddindex=1}^\nDoms \discRestriction\transpose \discPU \discBilForm_\ddindex^{-1}\discRestriction \right),
\end{equation*}
where $\mat{I}$ is the identity matrix and $\discPU$ is the matrix representation of the (discretized) multiplication by the partition of unity function $\pu_\ddindex$ from \cref{eq:pu}. With the above notation, the discrete MS-GFEM iteration for some initial guess $\vec{\fwp}^{(0)}$ is defined via
\begin{equation}\label{eq:msgfem:iteration:discrete}
  \vec{\fwp}^{(k+1)} := \vec{\fwp}^{(k)} + \discPreconditioner\discBilForm(\vec{\fwp}^e - \vec{\fwp}^{(k)}) = \vec{\fwp}^{(k)} + \discPreconditioner(\vec{\rhsF} - \discBilForm\vec{\fwp}^{(k)}), \quad k=0,1,\ldots,
\end{equation}
where $\vec{\rhsF}$ is the discretized right-hand side $\rhsF$ and $\vec{\fwp}^e$ is the exact solution of $\discBilForm\vec{\fwp}^e = \vec{\rhsF}$, the discretized PDE \cref{eq:weakFormGlobal}. In practice, this iteration is split into two steps:
\begin{subequations}
\label{eq:discrete:RAS}
\begin{alignat}{2}
  \label{eq:discrete:RAS1}
  &\vec{\fwp}^{(k+\frac{1}{2})} &&:= \vec{\fwp}^{(k)} + \sum_{\ddindex=1}^\nDoms \discRestriction\transpose \discPU \discBilForm_\ddindex^{-1}\discRestriction\left(\vec{\rhsF} - \discBilForm\vec{\fwp}^{(k)}\right),\\
  \label{eq:discrete:RAS2}
  &\vec{\fwp}^{(k+1)} &&\phantom{:}= \vec{\fwp}^{(k+\frac{1}{2})} + \discCoarseRestriction\transpose\discBilForm_S^{-1}\discCoarseRestriction\left(\vec{\rhsF} - \discBilForm\vec{\fwp}^{(k+\frac{1}{2})}\right),
\end{alignat}
\end{subequations}
i.e., a restricted additive Schwarz (RAS) method on the fine level, \cref{eq:discrete:RAS1}, and a subsequent multiplicative incorporation of a coarse space correction, \cref{eq:discrete:RAS2}. We thus refer to $\discPreconditioner$ as a two-level hybrid restricted additive Schwarz preconditioner. Note that incorporating the coarse space multiplicatively, which was motivated by \cref{eq:harmonicError}, provides a direct link between the MS-GFEM theory and the preconditioner, and is thus crucial for robustness.

Instead of performing the preconditioned Richardson-type iteration \cref{eq:msgfem:iteration:discrete}, we accelerate the solution of the (left-)preconditioned system $\discPreconditioner\discBilForm\vec{\fwp}^e = \discPreconditioner\vec{\rhsF}$ with a Krylov subspace method, which reduces the number of iterations at a similar computational cost per iteration. Since the preconditioner $\discPreconditioner$ is not symmetric, we choose GMRES for all simulations in \Cref{sec:experiments}.

\begin{remark}
An analysis of the discrete multiscale method (not the iterative version) for con\-vec\-tion-diffusion problems using conforming finite-element discretizations in a slightly different setting can be found in \cite[Section 6.2]{Ma2025}. 
In \cite{Strehlow2024,Ma2025Helmholtz}, an MS-GFEM-based preconditioner has been applied to elliptic and Helmholtz problems. Using conforming finite-element discretizations, the authors prove convergence of the discrete MS-GFEM iteration as well as of the preconditioned GMRES iteration with convergence rates depending on $\sqrt{\evalMax}$. However, our main goal is to demonstrate the practical use and flexibility of the method. Thus, in \Cref{sec:experiments} we carry out thorough numerical experiments on structured and unstructured grids with discontinuous Galerkin, finite volume, and conforming finite element discretizations. Since the discrete convergence theory depends on the specific discretization chosen, it is beyond the scope of this work to carry the analysis out for the different cases considered here. However, the continuous convergence result for the special case $\divergence \velField \equiv 0$ in \cref{th:msgfem:iteration} suggests that this should carry over to the discrete setting as in the symmetric, conforming cases \cite{Strehlow2024,Ma2025Helmholtz}. 
\end{remark}

\subsection{Discontinuous Galerkin discretization}
\label{sec:DG}

For most of the numerical simulations, we employ a weighted symmetric interior penalty discontinuous Galerkin (DG) discretization with upwinding 
(see e.g.\ \cite[Section 4.6]{DiPietroErn2011}). Hence, we outline the method below for completeness.

First, we introduce some notation. We collect all interior faces $\face = \partial \element^-(\face) \cap \partial \element^+(\face)$, where $\element^-(\face)$, $\element^+(\face)\in\mesh$, in the set $\mathcal{F}_h^I$. Similarly, all boundary faces $\face = \partial \element^-(\face) \cap \boundary$, where $\element^-(\face)\in\mesh$, are collected in the set $\mathcal{F}_h^{\partial\Omega}$, which we split into faces on the Dirichlet boundary, $\mathcal{F}_h^{D}\subseteq \mathcal{F}_h^{\partial\Omega}$, and faces on the outflow boundary, $\mathcal{F}_h^{O} \subseteq \mathcal{F}_h^{\partial\Omega}$. For each $\face\in\mathcal{F}_h^I$, we define the unit normal vector $\bn_\face$ to be oriented from $\element^{-}(\face)$ to $\element^+(\face)$. The unit normal of a boundary face $\face\in\mathcal{F}_h^{\partial\Omega}$ coincides with the outer unit normal to $\boundary$. Lastly, we define the positive part $(\velField\cdot\bn)^{\oplus} := \max(\velField\cdot\bn, 0)$ and negative part $(\velField\cdot\bn)^{\ominus} := \max(-\velField\cdot\bn, 0)$ of the normal flow.
The DG finite element space of degree $p$ on the mesh $\mesh$ is
\begin{equation*}
V_h^\text{DG} = \big\{ v\in L^2(\Omega) : v|_\element \in\mathbb{P}_p \; \forall \, \element \in \mesh\big\},
\end{equation*}
where $\mathbb{P}_p$ is either the set of polynomials of total degree $p$ for simplices
or the set of polynomials of maximum degree $p$ for cuboid elements. 
A function $v\in V_h^\text{DG}$ is two-valued on an interior face $\face\in\mathcal{F}_h^I$; we denote by $v^-$ and $v^+$ the restriction to $\face$ from $\element^-(\face)$ and $\element^+(\face)$, respectively.
For any point $\x\in \face \in \mathcal{F}_h^I$, we define the jump 
\begin{equation*}
  \llbracket v \rrbracket (\x) = v^-(\x)-v^+(\x)
\end{equation*}
and the weighted average
\begin{equation*}
  \{ v \}_w (\x) = w^- v^-(\x) + w^+ v^+(\x)
\end{equation*}
for some weights $w^- + w^+ = 1$, $w^\pm \geq 0$. If $w^- = w^+ = 1/2$, we simply write $\{v\}$.
A particular choice of weights depending on the
diffusion tensor $\diffTensor$ has been introduced in \cite{ErnStephansenZunino2008}. Assuming that
$\diffTensor^\pm$ is constant on $\element^\pm(\face)$, they set
$w^- = \delta_{\diffTensor\bn}^+/(\delta_{\diffTensor\bn}^- + \delta_{\diffTensor\bn}^+)$ and
$w^+ = \delta_{\diffTensor\bn}^-/(\delta_{\diffTensor\bn}^- + \delta_{\diffTensor\bn}^+)$ for
$\delta_{\diffTensor\bn}^{\pm} =\bn_\face^T \diffTensor^\pm \bn_\face$.

The DG method for numerically solving \cref{eq:fwp} now reads:
\begin{equation*}
\text{Find $u_h \in V_h^\text{DG}$ such that } 
  b_h^\text{DG}(u_h,v) =  l_h^\text{DG}(v)
  \text{ for all $v \in V_h^\text{DG}$,}
\end{equation*}
where
\begin{equation*}
  l_h^\text{DG}(v) = \sum_{\element\in\mesh} \int_\element fv \dx 
  - \sum_{\face\in\mathcal{F}_h^{O}} \int_\face \bcO v \ds 
  - \sum_{\face\in\mathcal{F}_h^D} \int_\face \bcD \Bigl[  \diffTensor\nabla v \cdot \bn_\face
  - \sigma_\face v
  + (\velField\cdot\bn_\face)^{\ominus} v\Bigr] \ds
\end{equation*}
and
\begin{equation*}
b_h^\text{DG}(u,v) = 
\sum_{\element\in\mesh} b_\element(u,v) 
+ \sum_{\face\in\mathcal{F}_h^I}  b^I_\face(u,v)
+ \sum_{\face\in\mathcal{F}_h^{\boundary}} b^{\boundary}_\face(u,v)
+ \sum_{\face\in\mathcal{F}_h^D} b^D_\face(u,v)
\end{equation*}
with 
\begin{align*}
b_\element(u,v) &= \int_\element \left(\diffTensor\nabla u - u \velField \right) \cdot \nabla v\dx,\\
b^I_\face(u,v) &= \int_\face \left(\sigma_\face + \frac{1}{2} \abs{\velField \cdot \bn_\face} \right) \llbracket u \rrbracket \llbracket v \rrbracket \ds 
+ \int_\face (\velField \cdot \bn_\face) \{u\} \llbracket v \rrbracket \ds \\ 
&\phantom{={}} - \int_\face \{\diffTensor\nabla u\}_{w} \cdot \bn_\face \llbracket  v \rrbracket \ds
- \int_\face \{\diffTensor\nabla v\}_{w} \cdot \bn_\face \llbracket  u \rrbracket \ds, \notag \\
b^{\boundary}_\face(u,v) &= \int_{\face} (\velField \cdot \bn_\face)^{\oplus} u v \ds,\\
b^D_\face(u,v) &= \sigma_\face \int_\face uv \ds 
- \int_\face \diffTensor\nabla u \cdot \bn_\face v \ds   
- \int_\face \diffTensor\nabla v \cdot \bn_\face u \ds.
\end{align*}
For each face, $\sigma_\face>0$ defines a penalty parameter, which is chosen as 
\begin{equation*}
  \sigma_\face = \begin{cases}
  \alpha \dfrac{2\delta_{\diffTensor\bn}^-\delta_{\diffTensor\bn}^+}{\delta_{\diffTensor\bn}^- + \delta_{\diffTensor\bn}^+} p(p + \dimDomain -1)\dfrac{\abs{\face}}{\min(\abs{\element^-(\face)}, \abs{\element^+(\face)})} & \text{ if }\; \face \in \mathcal{F}_h^I,\\[1em]
  \alpha \delta_{\diffTensor\bn} p(p + \dimDomain -1)\dfrac{\abs{\face}}{\abs{\element^-(\face)}} & \text{ if }\; \face \in \mathcal{F}_h^{\boundary},
  \end{cases}
\end{equation*}
where $\abs{\face}$ denotes the volume of $\face$ and $\alpha > 0$ is some user-defined parameter. This choice of penalty parameter was used in the benchmark paper \cite{Bastian2011} and is based on a combination of different approaches: The harmonic average of diffusion values was introduced and analyzed in \cite{ErnStephansenZunino2008}, the dependence on the polynomial degree was investigated in \cite{EpshteynRiviere2007}, and the specific form of mesh-dependence is taken from \cite{HartmannHouston2008}. In all numerical experiments of \Cref{sec:experiments}, we choose $\alpha = 3$ and $p=1$.

\subsection{Implementation and hardware}
\label{sec:implementation}

The MS-GFEM-based two-level RAS preconditioner described in this paper has been implemented within the DUNE software framework\footnote{www.dune-project.org} \cite{Bastian2008,Bastian2021} in a partially multithreaded setting, where in particular the generalized eigenproblems and the subdomain solves are carried out in parallel. The eigenproblems are solved using the C++ library Spectra \cite{spectralib} in symmetric shift-invert mode. The local subdomain problems as well as the global coarse problem are solved using UMFPack \cite{umfpack}. Furthermore, we use ParMETIS \cite{metis} for graph partitioning. Runtimes reported below are in seconds and were obtained on a dual-socket AMD EPYC 7713 system with 64 cores on each socket.


\section{Numerical results}
\label{sec:experiments}

In this section, we demonstrate the numerical performance of our proposed method in terms of robustness with respect to model parameters and scalability in large-scale problems. As an extension, we demonstrate that our method also performs well in the limit when the diffusion tensor vanishes, i.e., we solve a hyperbolic transport equation. 

The fine mesh $\mesh$ is divided into $\nDoms$ non-overlapping subdomains, each of which is then extended by two layers of elements to create the overlapping subdomains $\ovlpdomain_\ddindex$, $\ddindex=1,\ldots,\nDoms$. These overlapping subdomains are then further extended by $\nOvs$ additional element layers to create the oversampling domains $\osdomain_\ddindex$, $\ddindex=1,\ldots,\nDoms$. If not stated otherwise, we use $\nOvs = 2$ oversampling layers. For all simulations, we choose a distance-based partition of unity as defined in \cite[p.\ 84]{ToselliWidlund2005}.

We employ different strategies for choosing the coarse space size: Either a fixed number $\locDim_\ddindex = \evalFixed \in \N$ of basis functions per subdomain is preselected, which yields control over the coarse space size, or we set a threshold $\evalMax > 0$ and choose all eigenfunctions corresponding to larger eigenvalues, yielding an adaptive selection criterion with control over the convergence rate. Note the slight abuse of notation since this actually means that $\evalMax$ from \cref{eq:evalMax} is smaller than the threshold we choose.

GMRES is restarted after $100$ iterations and stopped once
\begin{equation*}
\big\vert\discPreconditioner(\vec{\rhsF}-\discBilForm \vec{\fwp}^{(k)})\big\vert < \errTol \, \big\vert\discPreconditioner(\vec{\rhsF}-\discBilForm \vec{\fwp}^{(0)})\big\vert,
\end{equation*}
i.e., a (preconditioned) residual reduction of $\errTol > 0$ has been reached; we then report the number of GMRES iterations $\nIter = k$ needed. The particular choice of $\errTol$ is specified for each experiment individually.


\subsection{Robustness with respect to grid P\'{e}clet number and subdomains}
\label{sec:checkerboard}

\begin{figure}
\centering
\begin{tikzpicture}[baseline=(image.base)]%
  \node (image) at (0,0) {\includegraphics[width=.425\textwidth]{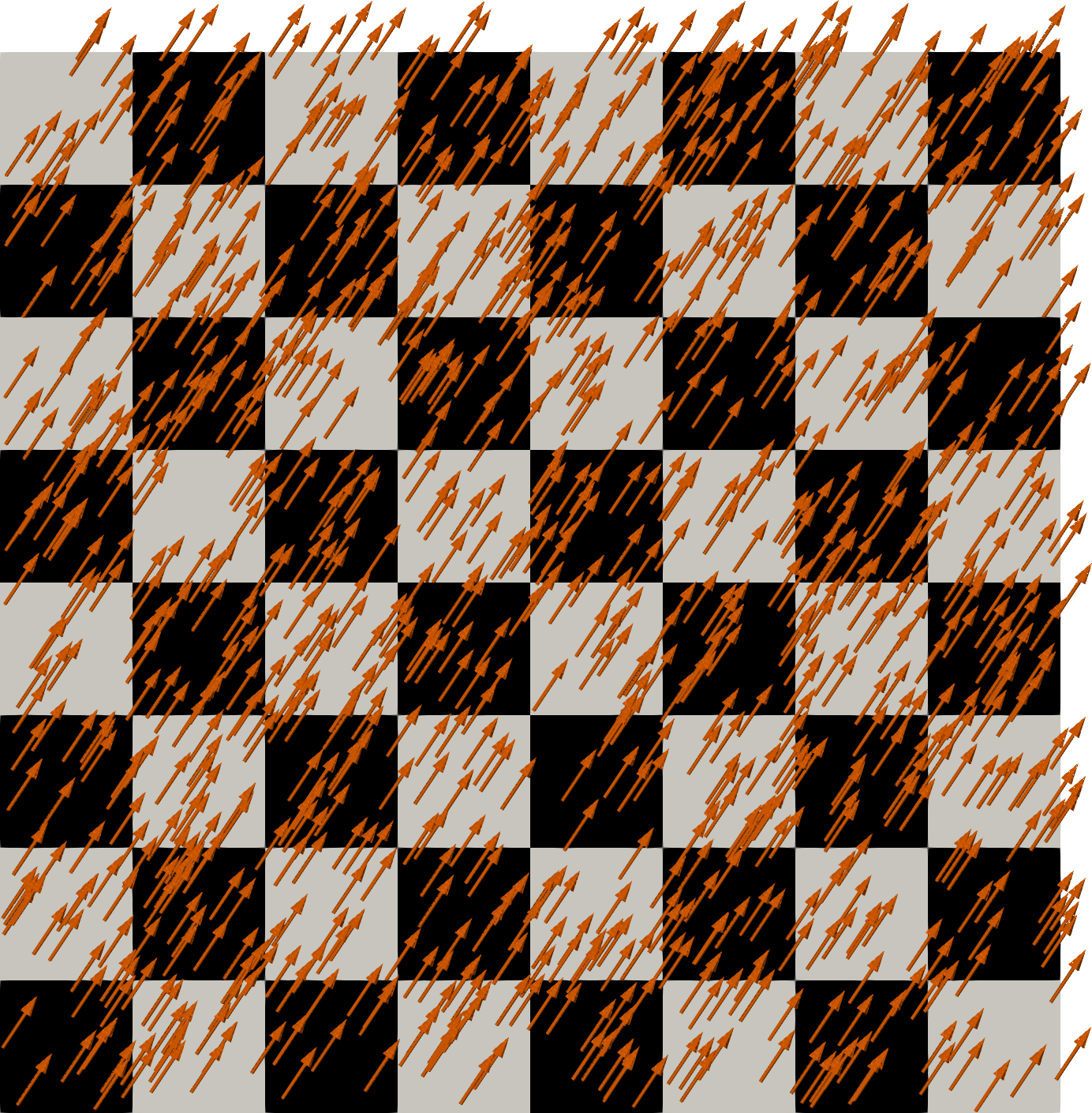}};%
\end{tikzpicture}%
\quad
\begin{tikzpicture}[baseline=(image.base)]%
  \node (image) at (0,0) {\includegraphics[width=.425\textwidth]{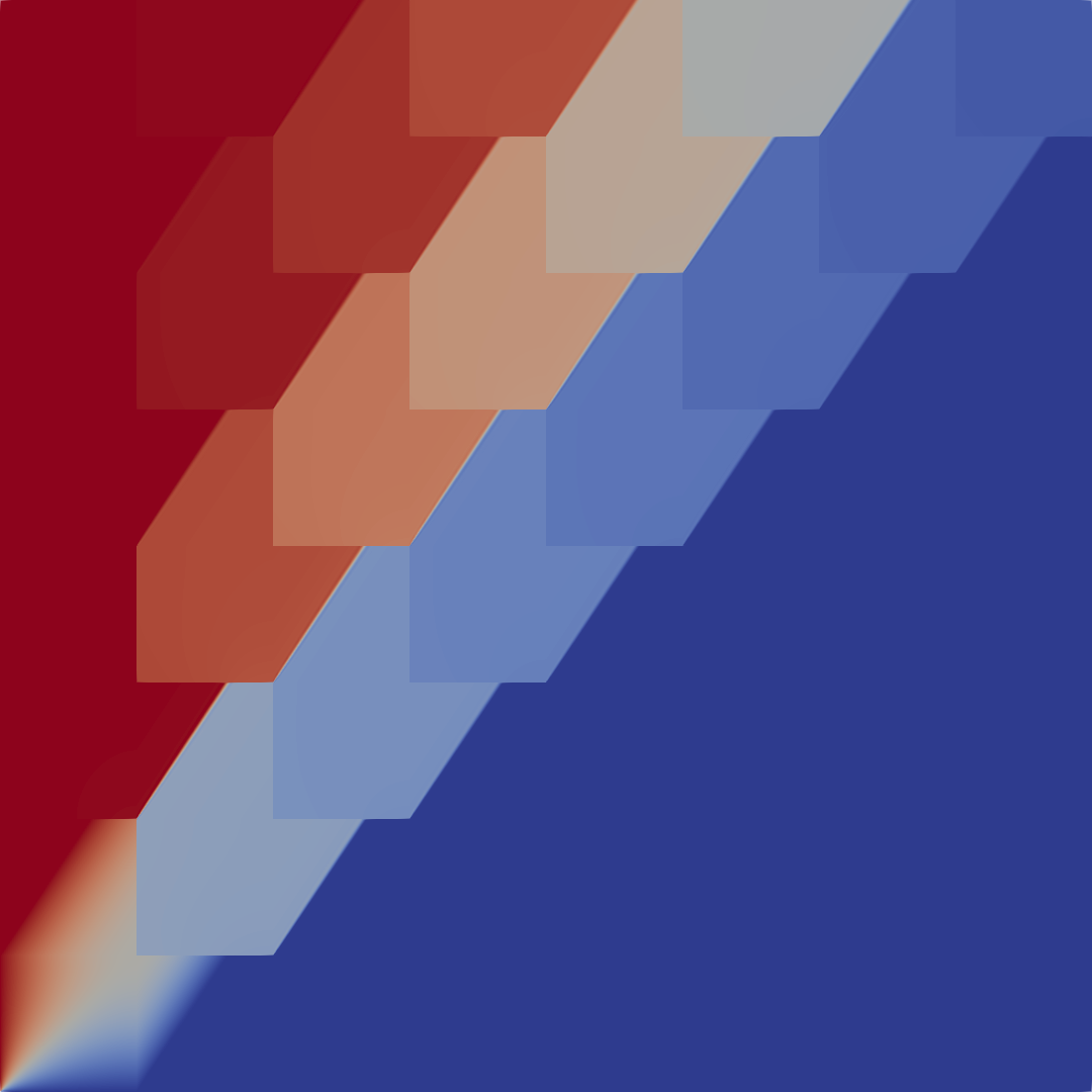}};%
  \node[anchor=center, yshift=4pt] at (image.north) {{{$\diffTensor \nabla \fwp \cdot \bn = 0$}}};%
  \node[anchor=center, yshift=-2pt] at (image.south) {{{$\fwp = 0$}}};%
  \node[anchor=center,rotate=90, yshift=-4pt] at (image.east) {{{$\diffTensor \nabla \fwp \cdot \bn = 0$}}};%
  \node[anchor=center,rotate=90, yshift=2pt] at (image.west) {{{$\fwp = 1$}}};%
\end{tikzpicture}%
\caption{Model setup for \Cref{sec:checkerboard}. Left: Diffusion coefficient and velocity field. The diffusion coefficient takes values of $\diffTensor \equiv \diffMax = 1$ and $\diffTensor \equiv \diffMin$ in the black and grey squares, respectively. Right: Corresponding PDE solution for $\diffMin = 10^{-6}$.}
\label{fig:checkerboard:model}
\end{figure}

We first consider a heterogeneous high-contrast model with constant velocity field on the domain $\domain = [0,1]^2$, which is discretized using a mesh of $1\,000 \times 1\,000$ equally-sized square elements. The diffusion coefficient is piecewise constant with 8 by 8 tiles in a checkerboard pattern alternating between the values $\diffMax = 1$ and $\diffMin \ll 1$, for which we consider different values in the simulations. The velocity field is defined as
\begin{equation}\label{eq:checkerboard:velField}
  \velField :\equiv \begin{pmatrix}2/3 \\ 1\end{pmatrix},
\end{equation}
and we choose $\bcO \equiv 0$ on the outflow boundary $\bdryO = (\{1\} \times [0,1]) \cup ([0,1] \times \{1\})$. On the Dirichlet part of the boundary, $\bdryD = \boundary \setminus \bdryO$, we set
\begin{equation}\label{eq:checkerboard:bc}
  \bcD(x_1,x_2) = \begin{cases}
  1 & \text{ if }\; x_1 = 0,\\
  0 & \text{ if }\; x_2 = 0,
  \end{cases}
\end{equation}
and there are no source terms in the domain, i.e., $\rhssource \equiv 0$. The model parameters and the DG solution for $\diffMin = 10^{-6}$ are visualized in \Cref{fig:checkerboard:model}. While the angle of the velocity field \cref{eq:checkerboard:velField} is arbitrary, we confirmed in simulations with a range of different angles between $0$ and $90$ degrees that the resulting GMRES iteration numbers and coarse space sizes are barely affected by the angle and the one we chose is thus representative.

\begin{figure}
\centering
\includegraphics[width=.49\textwidth]{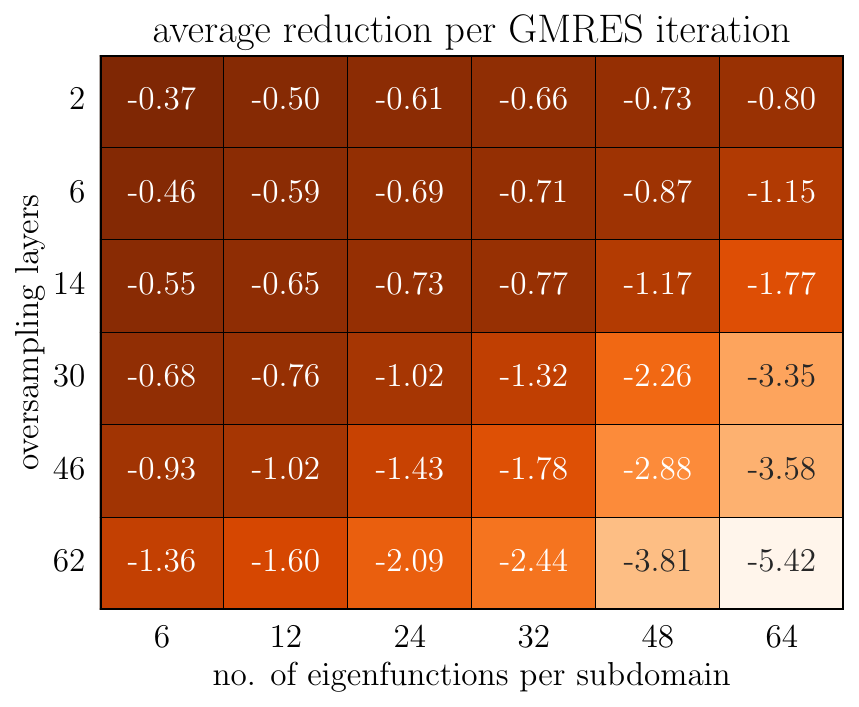}
\caption{Convergence test for the model of \Cref{sec:checkerboard} with $\diffMin = 10^{-6}$ using $100$ equally sized subdomains. The values indicate the decadic logarithm of the average residual reduction per GMRES iteration for different numbers $\evalFixed$ and $\nOvs$ of  eigenfunctions per subdomain and oversampling layers, respectively.}
\label{fig:checkerboard:convtest}
\end{figure}

In order to investigate the convergence rate of our preconditioner, we partition the domain into $10 \times 10$ subdomains of the same size and solve the system for different oversampling layers $\nOvs$ and coarse space sizes. For this convergence test, we choose a fixed number of eigenfunctions, $\evalFixed$, per subdomain. \Cref{fig:checkerboard:convtest} lists the decadic logarithm of the average residual reduction per GMRES iteration as a function of $\nOvs$ and $\evalFixed$, where we observe convergence rates of more than $10^{-5}$. As expected, the method converges faster with increasing $\nOvs$ and $\evalFixed$. Note that more oversampling layers result in more expensive local solves as well as eigenproblems, whereas increasing the number of eigenfunctions yields a larger coarse space and thus more costly coarse solves. It is dependent on the specific problem which values are optimal to solve the system in minimal time. However, we also obtain a robust preconditioner with small oversampling and coarse space dimension as we show next. To balance oversampling and coarse space size, we use the adaptive criterion of choosing $\evalMax$ for the remainder of this section. 

\renewcommand{\topfraction}{0.9}	

\begin{figure}
\centering
\includegraphics[width=.49\textwidth]{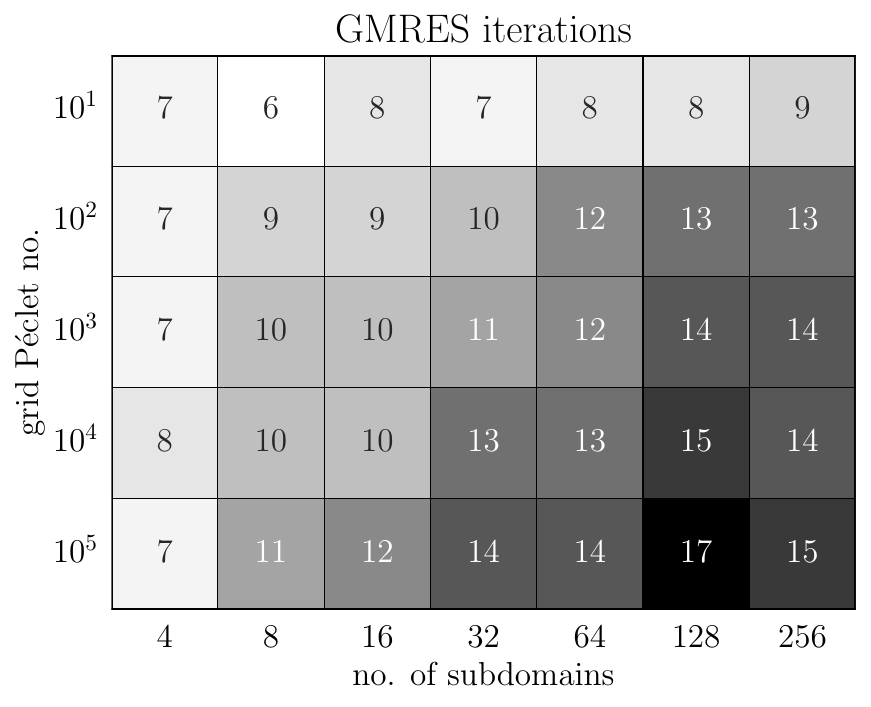}
\includegraphics[width=.49\textwidth]{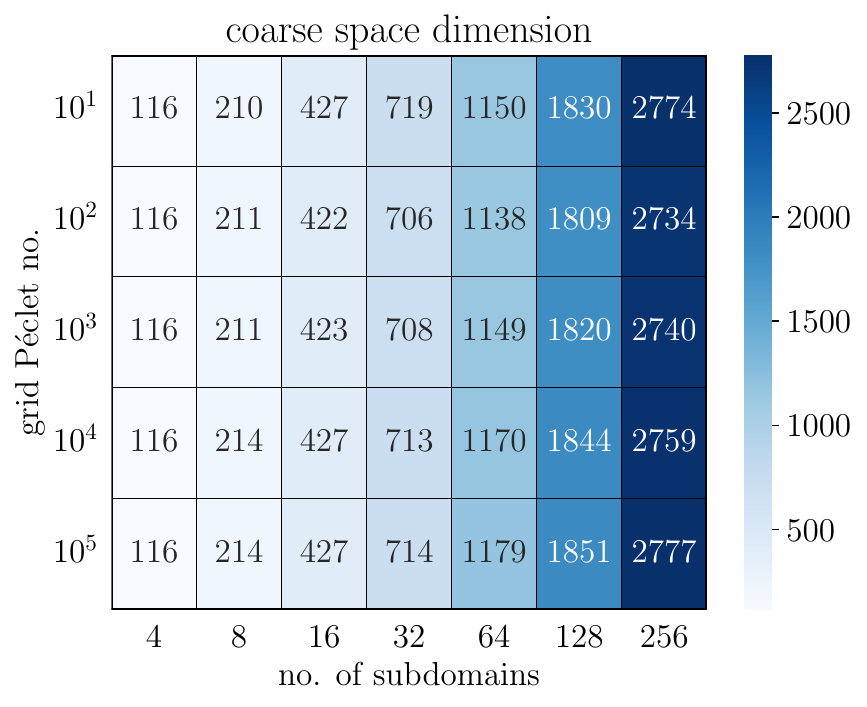}
\caption{Robustness with respect to grid P\'{e}clet number $\peclet_h$ and number of subdomains for the model of \Cref{sec:checkerboard}. 
We vary the grid P\'{e}clet number by choosing different values for $\diffMin$. An adaptive choice of eigenfunctions with $\evalMax = 2$ is used; all simulations are run until a residual reduction of $\errTol = 10^{-6}$ is achieved within GMRES. Left: Number of GMRES iterations $\nIter$. Right: Coarse space dimensions $\locDim$.}
\label{fig:checkerboard:robustness}
\end{figure}

A crucial property of our method is its robustness with respect to changes in the grid P\'{e}clet number and the number of subdomains. To demonstrate this, we use ParMETIS to partition the domain into $\nDoms \in \{4,8,16,32,64,128,256\}$ subdomains and vary the lower value of the diffusion coefficient $\diffMin$ to obtain grid P\'{e}clet numbers between $10^1$ and  $10^5$. We fix a small oversampling size of $\nOvs = 2$ and use the threshold $\evalMax = 2$ for the adaptive criterion. \Cref{fig:checkerboard:robustness} reports the coarse space dimension (right) and the number of GMRES iterations needed to achieve a residual reduction of $\errTol = 10^{-6}$ (left). The iteration numbers stay low for all different parameter combinations, underlining the robustness of the preconditioner even in the case with little oversampling and low coarse space dimension. Furthermore, the coarse space dimension is almost invariant with respect to changes in the grid P\'{e}clet number, and even for $256$ subdomains the coarse space dimension does not grow larger than $0.07\%$ of the DOFs of the global fine-scale discretization. Notably, the number of eigenfunctions per subdomain decreases the more subdomains are used, resulting in cheaper eigenvalue problems. This is in accordance with \cref{th:msgfem:evaldecay}: Since we keep the overlap and oversampling fixed at two layers each, the relative oversampling size becomes larger the more subdomains we use.

To demonstrate the flexibility of our preconditioner, we run simulations based on different fine-scale discretizations. We compare the DG discretization from \Cref{sec:DG} both with a cell-centered finite-volume (CCFV) scheme (e.g.\ \cite{LazarovMishevVassilevski1996}) and a conforming finite-element discretization with continuous piecewise linear basis functions (CG-$\mathbb{Q}1$, e.g.\ \cite{ErnGuermond2004}). As before, all schemes are based on a uniform mesh with $1\,000$ by $1\,000$ squares. For the comparison, we set $\diffMin = 10^{-5}$, i.e., $\peclet_h = 100$, and run GMRES until a residual reduction of $\errTol = 10^{-6}$ is achieved.

\Cref{tab:checkerboard:discComp} lists the iteration numbers as well as relative coarse space sizes for all three discretizations and different numbers of subdomains. While the CCFV scheme exhibits similar behavior as the DG method for $\evalMax = 2$, the iteration number seems to increase in this case when using CG-$\mathbb{Q}$1 elements, which is likely caused by the unstable conforming discretization method. Choosing a smaller threshold $\evalMax$ improves iteration numbers considerably, however, at the cost of a much larger coarse space. Thus, robustness can be enforced even for the unstable CG-$\mathbb{Q}$1 discretization. Note that, in general, the discrete theory depends on the specific discretization chosen and thus, we do not expect the same iteration count for a fixed threshold with different discretizations.

\begin{table}
\centering
\begin{tabular}{r | cc cc cc cc}
\toprule
 & \multicolumn{2}{c}{DG} & \multicolumn{2}{c}{CCFV} & \multicolumn{2}{c}{CG-$\mathbb{Q}1$} & \multicolumn{2}{c}{CG-$\mathbb{Q}1$}\\
  & \multicolumn{2}{c}{$\evalMax = 2$} & \multicolumn{2}{c}{$\evalMax = 2$} & \multicolumn{2}{c}{$\evalMax = 2$} & \multicolumn{2}{c}{$\evalMax = 0.5$}\\
$\nDoms$ & $\locDim$ & $\nIter$ & $\locDim$ & $\nIter$ & $\locDim$ & $\nIter$ & $\locDim$ & $\nIter$\\
\midrule
   4  & $0.003\%$ &  7 & $0.006\%$ &  9 & $0.006\%$ &  36 & $0.097\%$ & 11 \\
   8  & $0.005\%$ &  9 & $0.013\%$ & 10 & $0.011\%$ &  48 & $0.147\%$ & 15 \\
  16  & $0.011\%$ &  9 & $0.021\%$ & 11 & $0.020\%$ &  39 & $0.224\%$ & 19 \\
  32  & $0.018\%$ & 10 & $0.031\%$ & 12 & $0.028\%$ &  54 & $0.409\%$ & 14 \\
  64  & $0.028\%$ & 12 & $0.048\%$ & 13 & $0.045\%$ &  63 & $0.545\%$ & 16 \\
 128  & $0.045\%$ & 13 & $0.075\%$ & 13 & $0.070\%$ &  92 & $0.810\%$ & 24 \\
 256  & $0.068\%$ & 13 & $0.114\%$ & 13 & $0.102\%$ & 100 & $1.202\%$ & 26 \\
\bottomrule
\end{tabular}
\caption{Performance comparison of the preconditioner for discontinuous Galerkin (DG), cell-centered finite volume (CCFV), and conforming $\mathbb{Q}1$ (CG-$\mathbb{Q}1$) discretizations, using the model from \Cref{sec:checkerboard} with $\diffMin = 10^{-5}$, i.e., grid P\'{e}clet number $\peclet_h = 100$. GMRES is run until a residual reduction of $\errTol = 10^{-6}$ is achieved. We list iteration counts $\nIter$ and coarse space sizes $\locDim$ (in percent of the degrees of freedom on the fine mesh) for different numbers of subdomains $\nDoms$.
\label{tab:checkerboard:discComp}}
\end{table}


\subsection{A large-scale problem with non-divergence-free velocity}
\label{sec:largescale}

Next, we consider a large-scale problem and demonstrate robustness even for extreme cases with more than $10^5$ subdomains. The domain $\domain = [0,1]^2$ is discretized by a uniform grid with $5\,120$ squares in each direction, resulting in a total number of $104\,857\,600$ degrees of freedom for the fine DG discretization. In this example, we define the piecewise constant diffusion coefficient in a checkerboard pattern with 32 by 32 equally-sized tiles similar to \Cref{sec:checkerboard}, alternating between the values $\diffMax = 1$ and $\diffMin = 10^{-6}$. We choose the rotating velocity field
\begin{equation*}
  \velField(x_1,x_2) := \begin{pmatrix}
  20(x_2 - 0.5)[1 - (x_1 - 0.5)^2] - 1.5(x_1 - 0.5)\\
  20(x_1 - 0.5)[1 - (x_2 - 0.5)^2] - 1.5(x_2 - 0.5)
  \end{pmatrix},
\end{equation*}
which has negative divergence and yields an indefinite problem with a grid P\'{e}clet number of approximately $3\,000$. We employ homogeneous Dirichlet boundary conditions everywhere on $\boundary$ and define a Gaussian-like source term
\begin{equation*}
  \rhssource(x_1,x_2) := \exp\left(-\frac{(x_1 - 0.25)^2 + (x_2 - 0.5)^2}{10^{-2}}\right).
\end{equation*}
In all simulations, we partition the domain in equally-sized subdomains and choose the number of basis functions per subdomain adaptively via different threshold parameters $\evalMax$. In addition, we also compare our method with another two-level hybrid RAS preconditioner, in which the coarse space is defined by multiplying the constant and linear functions with the partition of unity on each subdomain such that we obtain $\dimDomain + 1$ basis functions per subdomain. We refer to this as the partition of unity (PoU) coarse space (see e.g.\ \cite[Section 4.2]{Sarkis2003,DoleanJolivetNataf2015}). In all cases, GMRES is run until a residual reduction of $\errTol = 10^{-8}$ or a maximum number of $1\,000$ iterations is achieved.

We observe robustness of our method for very high numbers of subdomains even when using lower-dimensional coarse spaces (\Cref{tab:largescale}). By choosing smaller thresholds $\evalMax$, the number of GMRES iterations can be kept smaller at the cost of a higher-dimensional coarse space. On the other hand, the preconditioner with the PoU coarse space converges very slowly. When the method is terminated after $1\,000$ iterations, a residual reduction of only $10^{-2}$ is obtained in all cases. 

\begin{table}
\centering
\begin{tabular}{r | cc cc cc | cc }
\toprule
& \multicolumn{6}{|c|}{MS-GFEM} & \multicolumn{2}{c}{PoU} \\
\midrule
 & \multicolumn{2}{c}{$\evalMax=2$} & \multicolumn{2}{c}{$\evalMax=1$} & \multicolumn{2}{c|}{$\evalMax=0.5$} &  &  \\
$\nDoms$ & $\locDim$ & $\nIter$ & $\locDim$ & $\nIter$ & $\locDim$ & $\nIter$ & $\locDim$ & $\nIter$\\
\midrule
 $15\,625$  & $0.12\%$ & 23 & $0.19\%$ & 19 & $0.41\%$ & 10 & $0.04\%$ & ($>1\,000$) \\
 $30\,625$  & $0.18\%$ & 28 & $0.28\%$ & 23 & $0.56\%$ & 10 & $0.09\%$ & ($>1\,000$) \\
 $62\,500$  & $0.25\%$ & 35 & $0.41\%$ & 30 & $0.80\%$ & 13 & $0.18\%$ & ($>1\,000$) \\
$122\,500$  & $0.39\%$ & 42 & $0.55\%$ & 36 & $1.08\%$ & 17 & $0.35\%$ & ($>1\,000$) \\
\bottomrule
\end{tabular}
\caption{Indefinite large-scale example from \Cref{sec:largescale} with $104\,857\,600$ degrees of freedom on the fine mesh and a grid P\'{e}clet number of approximately $\peclet_h \approx 3\,000$. For different numbers of subdomains $\nDoms$, we list again iteration counts $\nIter$ and coarse space sizes $\locDim$ in percent of the degrees of freedom on the fine mesh.}
\label{tab:largescale}
\end{table}


\subsection{3D system on unstructured grid}
\label{sec:3d}

\begin{figure}
\begin{center}
\begin{tikzpicture}[baseline=(image.base)]%
  \node (image) at (0,0) {\includegraphics[width=.475\columnwidth]{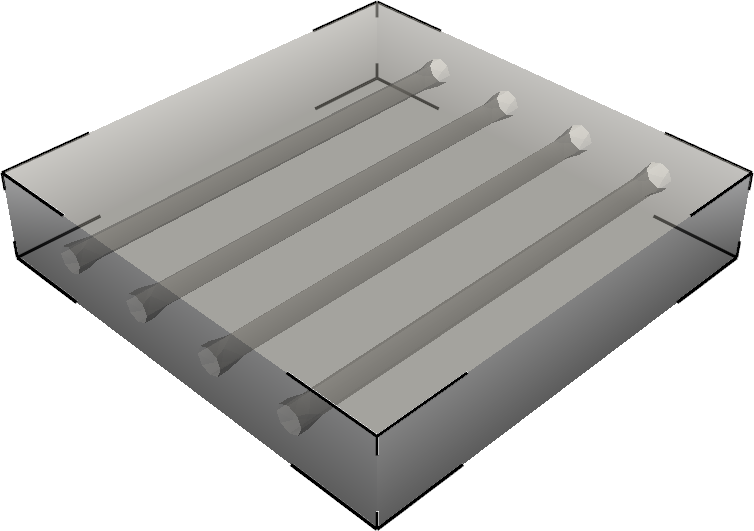}};%
  \node at (-3.0,2.25) {{(a)}};%
\end{tikzpicture}%
\begin{tikzpicture}[baseline=(image.base)]%
  \node (image) at (0,0) {\includegraphics[width=.475\columnwidth]{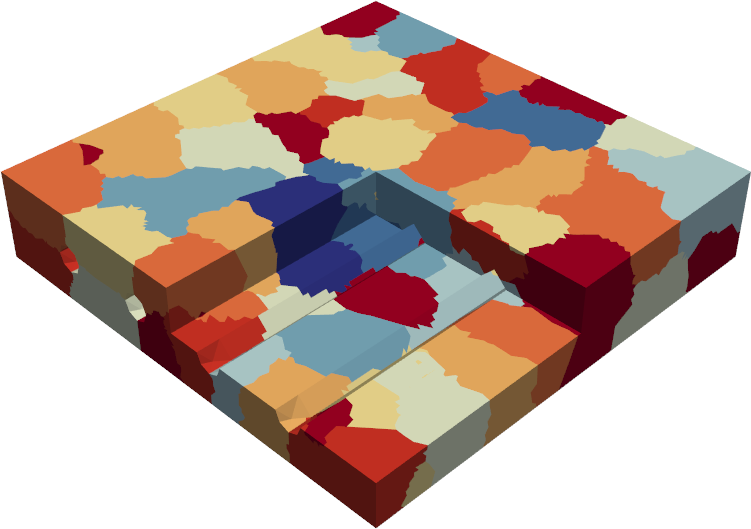}};%
  \node at (-3.0,2.25) {{(b)}};%
\end{tikzpicture}%
\linebreak
\begin{tikzpicture}[baseline=(image.base)]%
  \node (image) at (0,0) {\includegraphics[width=.925\columnwidth]{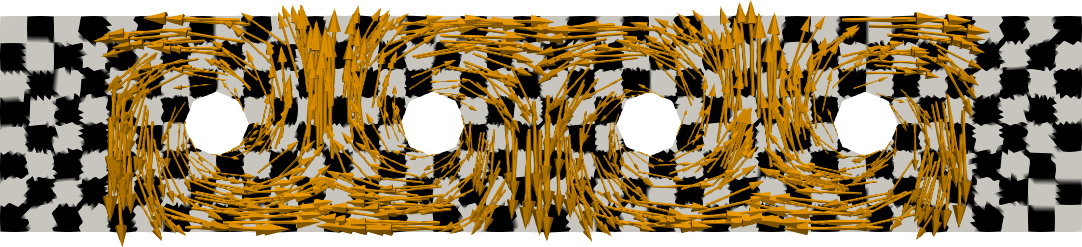}};%
  \node[anchor=center, xshift=-6pt] at (image.west) {{(c)}};%
\end{tikzpicture}%
\end{center}
\caption{Model setup for the 3D example from \Cref{sec:3d}. (a) Computational domain. (b) ParMETIS partitioning into 64 subdomains. (c) Vertical slice through the diffusion coefficient at $x_1=0$ and rotating velocity field.}
\label{fig:3d}
\end{figure}

To further underline the flexibility of the proposed method and to investigate its computational cost, we consider a 3D example with a complex domain (see \Cref{fig:3d}a). It consists of a rectangular prism $[0,10] \times [0,10] \times [-1,1]$ with four cylindrical cutouts along the $x_1$-axis with centers at $x_3 = 0$ and $x_2 \in \{2,4,6,8\}$ and is discretized using an unstructured tetrahedral mesh that contains $2\,016\,256$ degrees of freedom. The high-contrast diffusion coefficient alternates between the values $\diffTensor \equiv 1$ and $\diffTensor\equiv\diffMin$ in a 3D checkerboard pattern with $80 \times 40 \times 8$ rectangular prisms. We construct an artificial velocity field that rotates around the four pipes (see \Cref{fig:3d}c), choose $\rhssource \equiv 1$, and employ homogeneous Dirichlet boundary conditions everywhere on the boundary. We perform experiments with $\diffMin \in \{10^{-3}, 10^{-4}, 10^{-5}, 10^{-6}\}$, resulting in grid P\'{e}clet numbers ranging roughly between $40$ and $40\,000$.

For comparison, we export the assembled global systems to PETSc \cite{petsc} and solve them using GMRES with two different algebraic multigrid preconditioners: PETSc's GAMG and BoomerAMG from the hypre library \cite{hypre,RugeStueben1987}. Due to the DG discretization on the tetrahedral mesh, we set the block size of the PETSc matrix to 4. GAMG and BoomerAMG are then run with 64 MPI processes and the default settings from PETSc; tuning the parameters, like smoother type, connection threshold, or coarsening strategy, did not lead to better results than the ones listed below. 

For both hybrid Schwarz preconditioners, we partition the domain into 64 subdomains using ParMETIS (\Cref{fig:3d}b). These subdomains are extended by two layers of overlap and one additional layer of oversampling for MS-GFEM, whereas we choose an overlap of three layers of elements for the PoU method to obtain comparable results. The PoU coarse space is built from the constant and linear functions, i.e., we choose four basis functions per subdomain leading to a coarse space dimension of 256. For MS-GFEM, we compute a fixed number of only two eigenfunctions per subdomain to build a 128-dimensional coarse space. In both cases, the setup of the preconditioner as well as the local solves are executed in parallel using 64 threads, i.e., one thread per subdomain. 

All simulations are run until a residual reduction of $\errTol = 10^{-10}$ is achieved or a maximum number of iterations is reached. We stop the multigrid methods after $10\,000$ iterations and the two-level RAS methods after $1\,000$, at which point the runtime is already significantly larger than for the methods that converge. Note that all results were obtained on the same machine with the same number of CPU cores, which did not exceed the cores available on a single socket. Hence, the differences due to threads versus MPI processes are expected to be negligible.

\Cref{tab:3d_comparison} lists the number of GMRES iterations as well as the runtime of the solvers using the four different preconditioners for different values of $\diffMin$. Although GAMG and BoomerAMG are faster than the other methods if they converge, both of them only converge for $\diffMin = 10^{-3}$. While the PoU method is a bit more robust and converges for $\diffMin \in \{10^{-3}, 10^{-4}\}$, the runtime in the latter case is faster for the MS-GFEM preconditioner. In fact, the iteration number, setup time, as well as the total time needed is independent of the model parameter $\diffMin$ for MS-GFEM, highlighting the flexibility and robustness of the method once again. 
\begin{table}
\centering
\begin{tabular}{r | cr | cr | cr | ccc }
\toprule
& \multicolumn{2}{c|}{GAMG} & \multicolumn{2}{c|}{BoomerAMG} & \multicolumn{2}{c|}{PoU} & \multicolumn{3}{c}{MS-GFEM} \\
$\diffMin$ & $\nIter$ & $t_\mathrm{total}$ & $\nIter$ & $t_\mathrm{total}$ & $\nIter$ & $t_\mathrm{total}$ & $\nIter$ & $t_\mathrm{setup}$ & $t_\mathrm{total}$\\
\midrule
 $10^{-3}$  & $2\,061$ & 52.7 & 165 & 5.7 & 179 &  70.9 & 13 & 78.1 & 81.6\\
 $10^{-4}$  &    -     &   -  &  -  &  -  & 365 & 134.3 & 13 & 77.1 & 80.7\\
 $10^{-5}$  &    -     &   -  &  -  &  -  & -   &   -   & 12 & 77.4 & 80.4\\
 $10^{-6}$  &    -     &   -  &  -  &  -  & -   &   -   & 13 & 77.6 & 81.5\\
\bottomrule
\end{tabular}
\caption{3D example from \Cref{sec:3d} using different values 
of $\diffMin$. We list GMRES iteration counts $\nIter$, as well as total runtimes $t_\mathrm{total}$ and setup times $t_\mathrm{setup}$ for different preconditioners. All times are given in seconds. Dashes indicate that in that case 
a predefined maximum number of iterations was reached.}
\label{tab:3d_comparison}
\end{table}

As a further reference, we compute the solution of the equation using the parallel sparse direct solver MUMPS \cite{mumps} through PETSc with 64 MPI processes. Even for this moderately-sized example, the runtime lies between $300$ and $360$ seconds depending on $\diffMin$ and thus significantly above the MS-GFEM runtime. The difference is expected to become even bigger for larger problems, besides the infeasible memory requirements related to using direct solvers.

The dominating cost of the MS-GFEM preconditioner is clearly the setup phase, which involves solving the generalized eigenproblems. Recently, a more efficient version of MS-GFEM has been introduced by formulating these eigenproblems not on the full oversampling domains but only on the overlaps \cite{Alber2025}. In their numerical examples, solving the eigenproblems was accelerated by a factor of up to 8, significantly reducing the total cost of the method. While their theory is formulated for elliptic problems, it is a promising direction to investigate the performance of their method for non-elliptic PDEs.


\subsection{Performance near the hyperbolic limit}
\label{sec:transport}

\begin{figure}
\centering
\includegraphics[width=.49\textwidth]{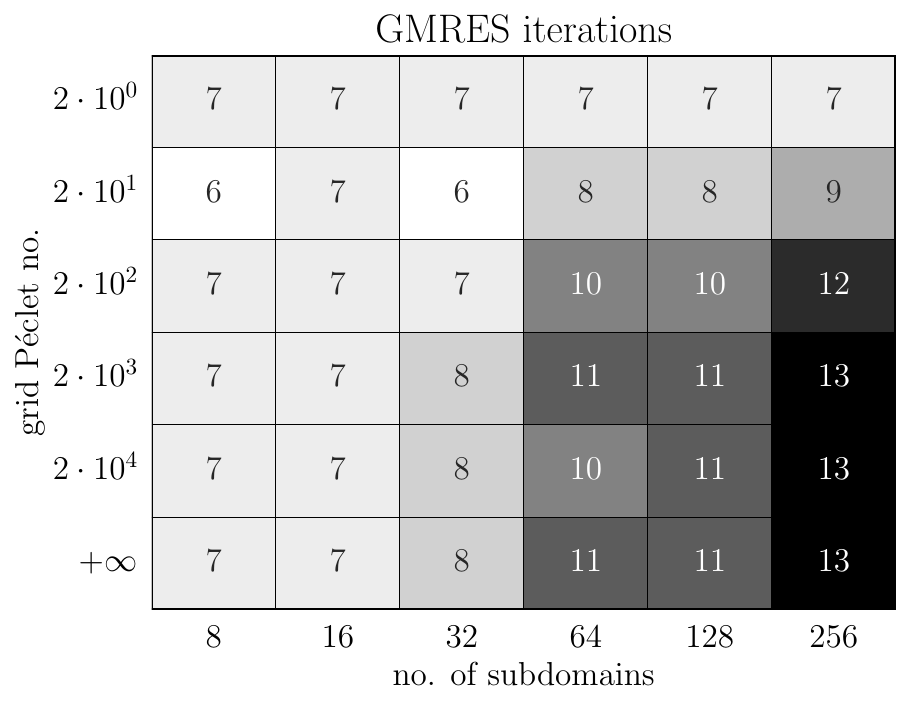}
\includegraphics[width=.49\textwidth]{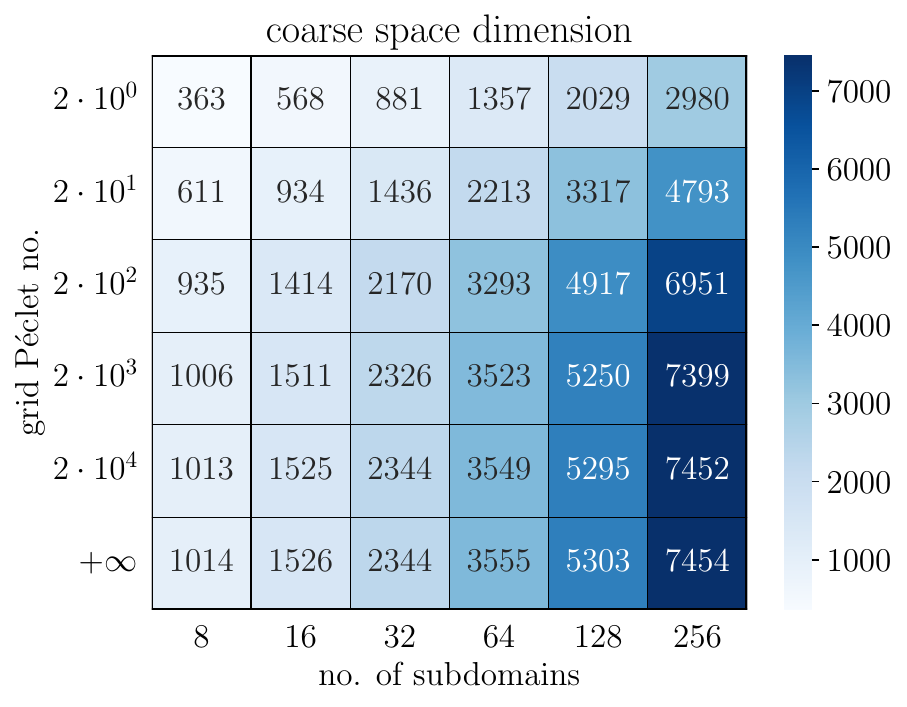}
\caption{Robustness in the pure-transport limit for the model of \Cref{sec:transport}, using constant diffusion and the 
velocity field from \Cref{sec:checkerboard}. In all simulations, $\evalMax = 0.5$ and GMRES is stopped once a residual reduction of $\errTol = 10^{-10}$ is achieved. Different grid P\'{e}clet numbers correspond to different values of the diffusion coefficient, where $+\infty$ indicates the choice $\diffTensor \equiv 0$. Left: Number of GMRES iterations $\nIter$. Right: Coarse space dimensions $\locDim$.}
\label{fig:transport}
\end{figure}

As an extension of our framework, we investigate how our method performs in the limit of a formally infinite P\'{e}clet number, i.e., we solve a first-order hyperbolic transport equation on the domain $\domain = [0,1]^2$, which is discretized by a uniform grid with $600$ squares in each direction. As in \Cref{sec:checkerboard}, we consider no source terms, $\rhssource \equiv 0$, and choose the velocity field \cref{eq:checkerboard:velField} and boundary conditions \cref{eq:checkerboard:bc}. The only difference to the model there is a constant diffusion coefficient $\diffTensor \equiv \diffMin$, which we vary and eventually set to $0$. In the latter case of pure convection, we need to specify an inner product for the restriction operator \cref{eq:resOp}. Since, in this particular case, the $B$-harmonic space does not contain diffusion, and any constant diffusion value in the inner product cancels out in the eigenvalue problem \cref{eq:gevpSaddlePoint}, we simply choose $A \equiv 1$ in the inner product, i.e., $\bilIP{\fwp}{\fwpTest}{\ovlpdomain} = \int_\ovlpdomain \nabla \fwp \cdot \nabla \fwpTest \dx$.

The number of GMRES iterations and the coarse space dimensions with the choice $\evalMax=0.5$ are reported in \Cref{fig:transport}. Once again, the number of iterations is robust with respect to changes in the grid P\'{e}clet number, even in the limit $\peclet_h = +\infty$. For a fixed number of subdomains, the coarse space dimension only changes with increasing grid P\'{e}clet numbers for moderately convection-dominated problems. In particular, there is almost no difference between the cases $\peclet_h \approx 2 \cdot 10^4$ and $\peclet_h = +\infty$, suggesting that all important coarse-space features are already captured in the case of high but finite grid P\'{e}clet numbers. This is probably due to the numerical approximation error and it is expected that refining the grid will shift this saturation of the coarse space to lower values of $\diffTensor$. Likewise, we ran the same simulations with the CCFV discretization and observed similar behavior. However, the saturation of the coarse space occurred for lower grid P\'{e}clet numbers, which is likely due to numerical diffusion and the generally larger approximation error compared to the DG discretization.

\begin{figure}
\begin{center}
  \begin{tikzpicture}[baseline=(image.base)]%
  \node (image) at (0,0) {\includegraphics[width=.15\textwidth]{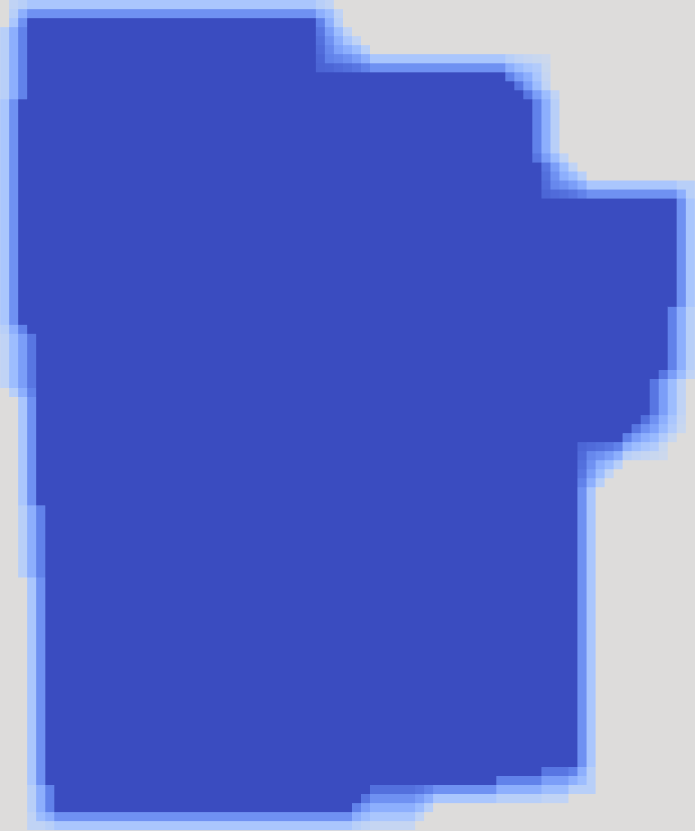}};%
  \node[anchor=center,rotate=90, yshift=2pt] at (image.west) {\footnotesize{{$\peclet_h = 0$}}};%
  \node[anchor=center, yshift=2pt] at (image.north) {\footnotesize{{$\evalindex=1$}}};%
  \end{tikzpicture}%
  \begin{tikzpicture}[baseline=(image.base)]%
  \node (image) at (0,0) {\includegraphics[width=.15\textwidth]{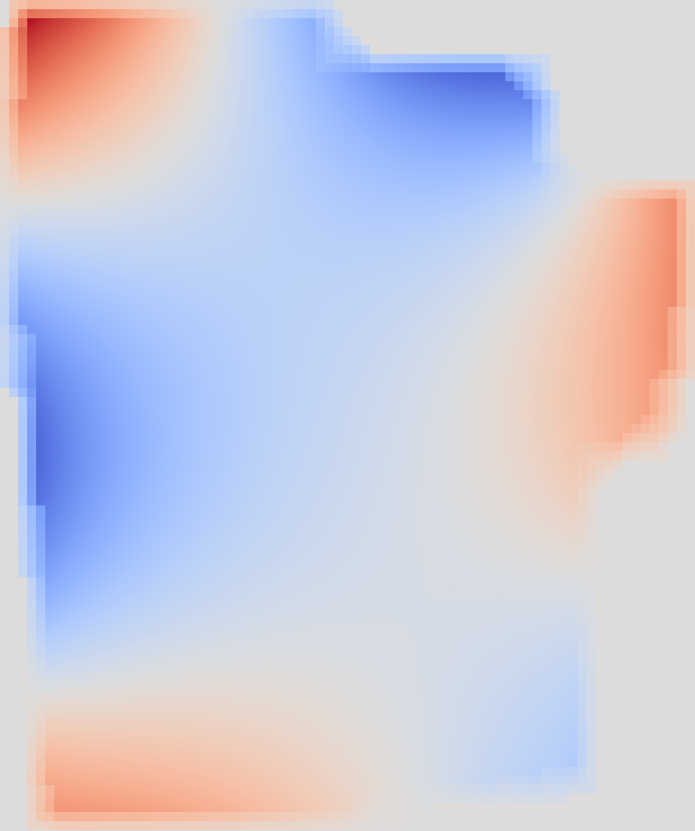}};%
  \node[anchor=center, yshift=2pt] at (image.north) {\footnotesize{{$\evalindex=5$}}};%
  \end{tikzpicture}%
  \begin{tikzpicture}[baseline=(image.base)]%
  \node (image) at (0,0) {\includegraphics[width=.15\textwidth]{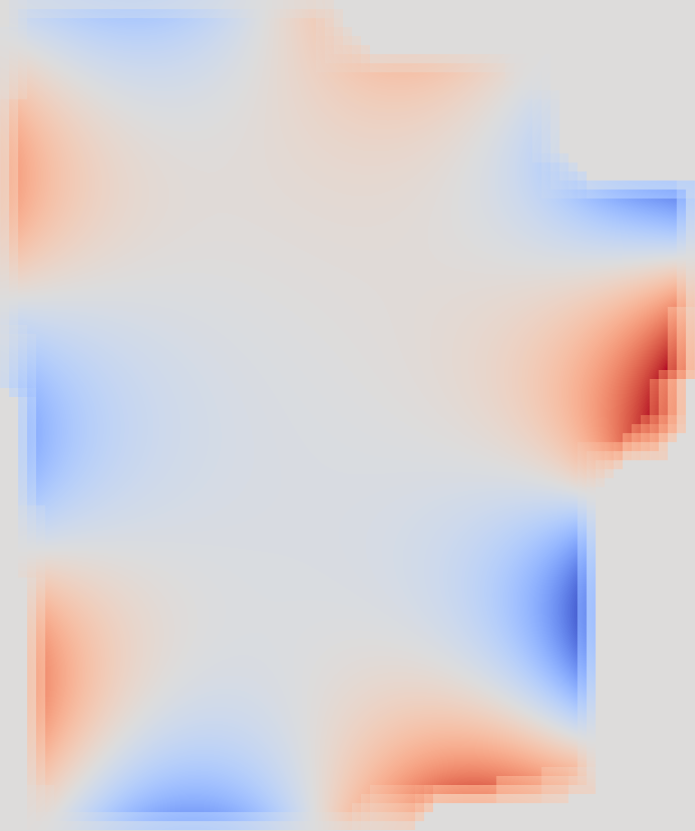}};%
  \node[anchor=center, yshift=2pt] at (image.north) {\footnotesize{{$\evalindex=10$}}};%
  \end{tikzpicture}%
  \begin{tikzpicture}[baseline=(image.base)]%
  \node (image) at (0,0) {\includegraphics[width=.15\textwidth]{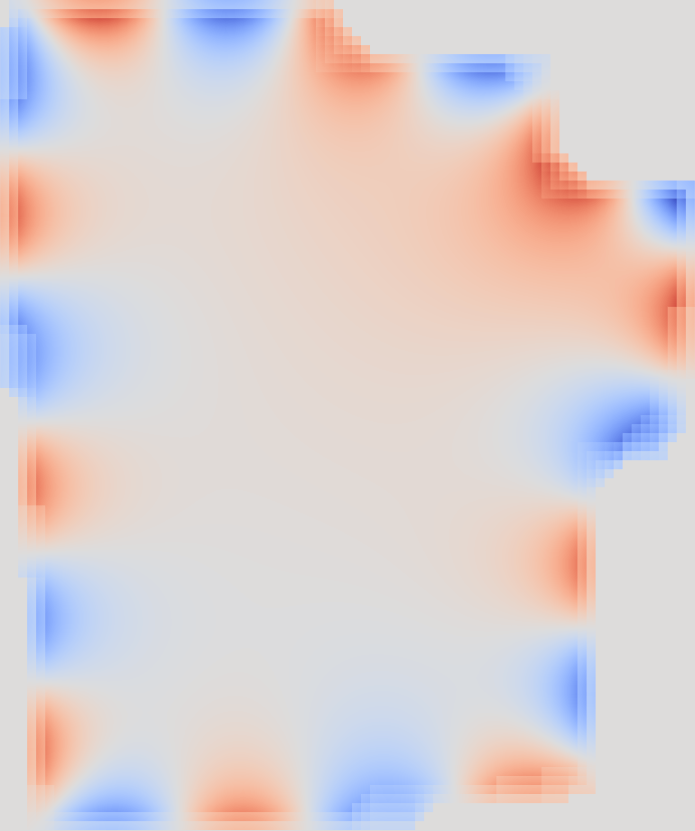}};%
  \node[anchor=center, yshift=2pt] at (image.north) {\footnotesize{{$\evalindex=20$}}};%
  \end{tikzpicture}%
  \begin{tikzpicture}[baseline=(image.base)]%
  \node (image) at (0,0) {\includegraphics[width=.15\textwidth]{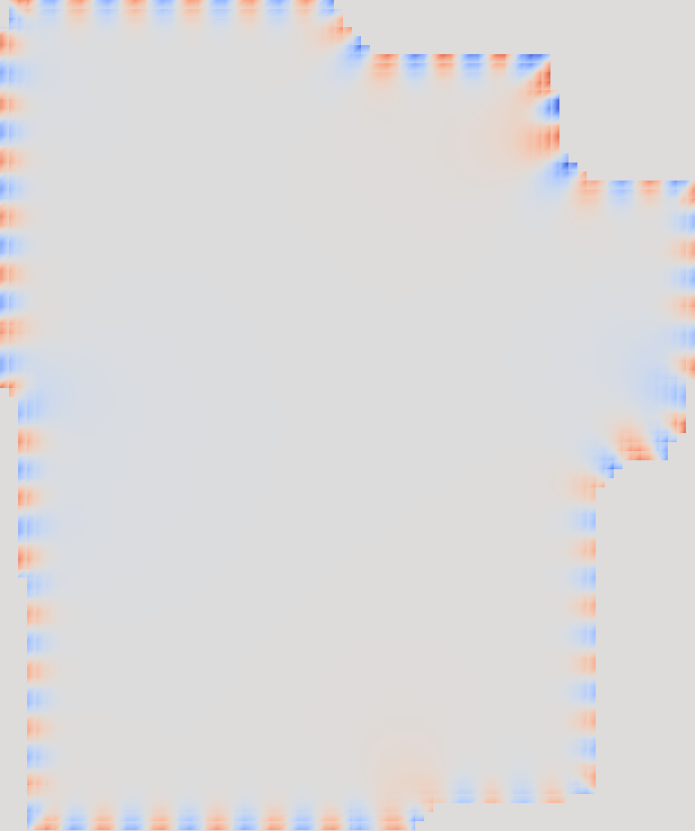}};%
  \node[anchor=center, yshift=2pt] at (image.north) {\footnotesize{{$\evalindex=100$}}};%
  \end{tikzpicture}%
  \linebreak
  \begin{tikzpicture}[baseline=(image.base)]%
  \node (image) at (0,0) {\includegraphics[width=.15\textwidth]{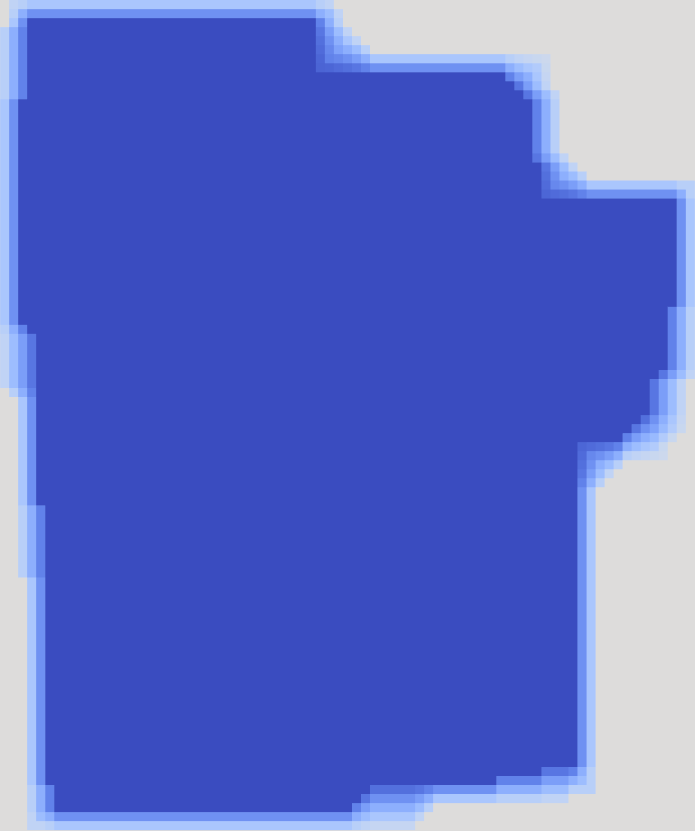}};%
  \node[anchor=center,rotate=90, yshift=2pt] at (image.west) {\footnotesize{{$\peclet_h = 2$}}};%
  \end{tikzpicture}%
  \begin{tikzpicture}[baseline=(image.base)]%
  \node (image) at (0,0) {\includegraphics[width=.15\textwidth]{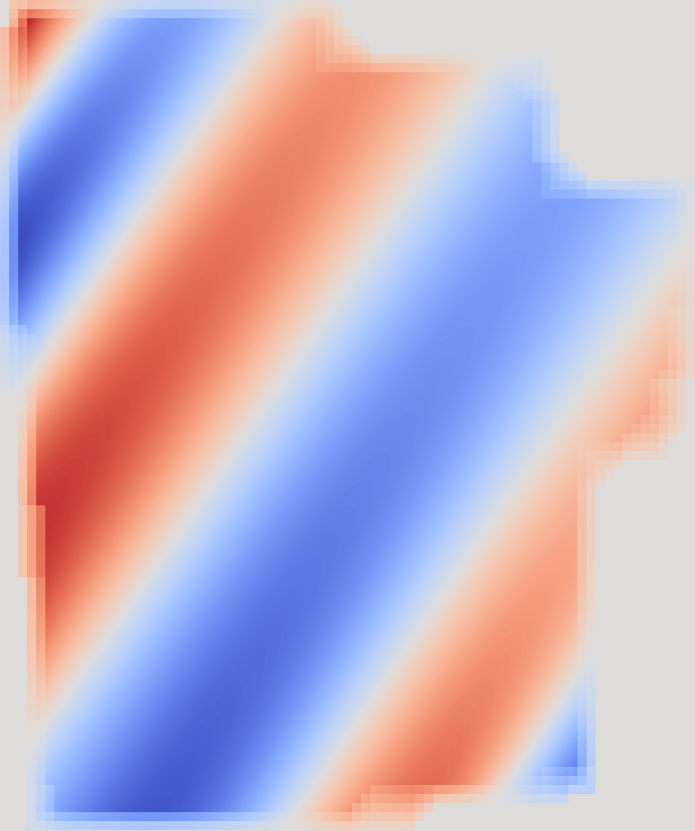}};%
  \end{tikzpicture}%
  \begin{tikzpicture}[baseline=(image.base)]%
  \node (image) at (0,0) {\includegraphics[width=.15\textwidth]{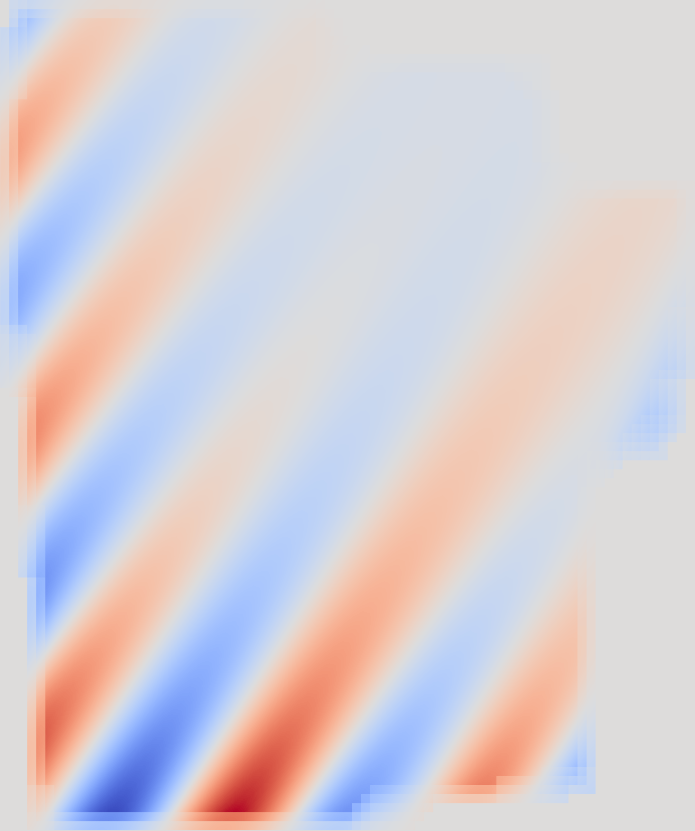}};%
  \end{tikzpicture}%
  \begin{tikzpicture}[baseline=(image.base)]%
  \node (image) at (0,0) {\includegraphics[width=.15\textwidth]{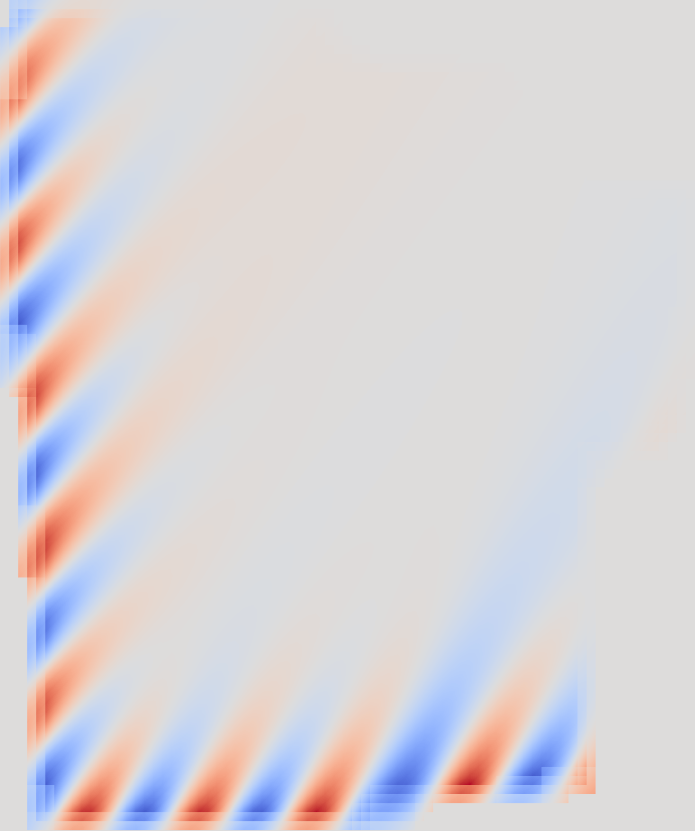}};%
  \end{tikzpicture}%
  \begin{tikzpicture}[baseline=(image.base)]%
  \node (image) at (0,0) {\includegraphics[width=.15\textwidth]{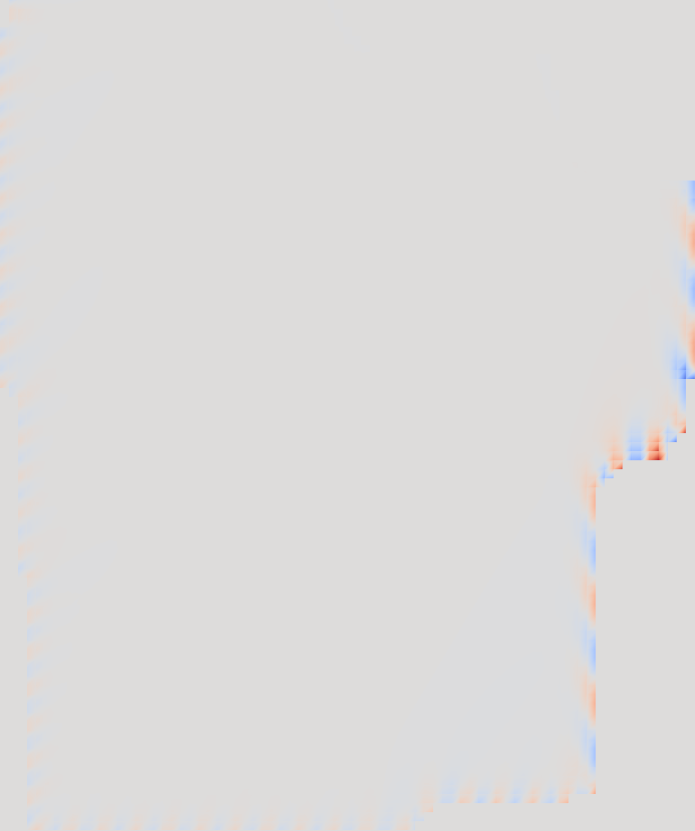}};%
  \end{tikzpicture}%
  \linebreak
  \begin{tikzpicture}[baseline=(image.base)]%
  \node (image) at (0,0) {\includegraphics[width=.15\textwidth]{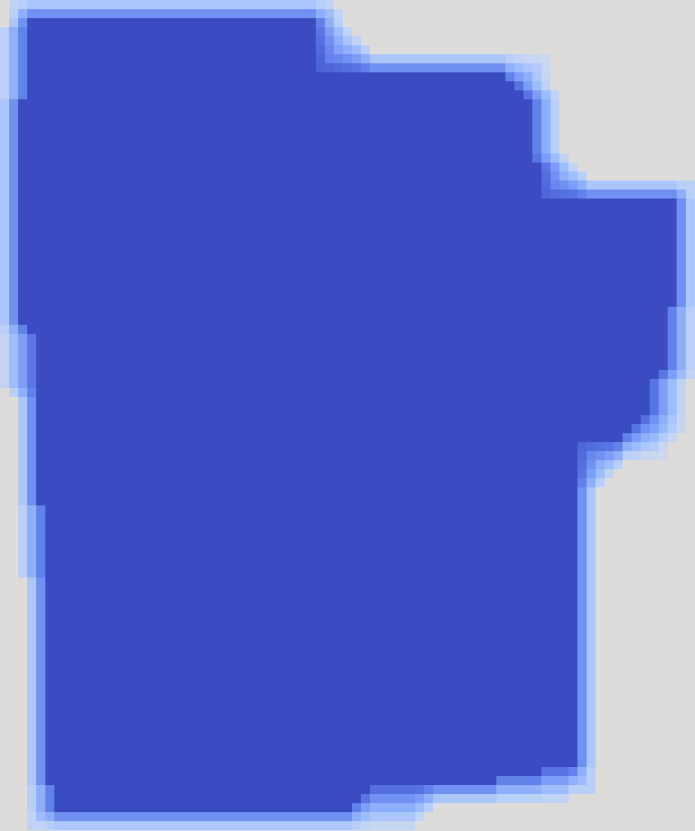}};%
  \node[anchor=center,rotate=90, yshift=2pt] at (image.west) {\footnotesize{{$\peclet_h = 2\,000$}}};%
  \end{tikzpicture}%
  \begin{tikzpicture}[baseline=(image.base)]%
  \node (image) at (0,0) {\includegraphics[width=.15\textwidth]{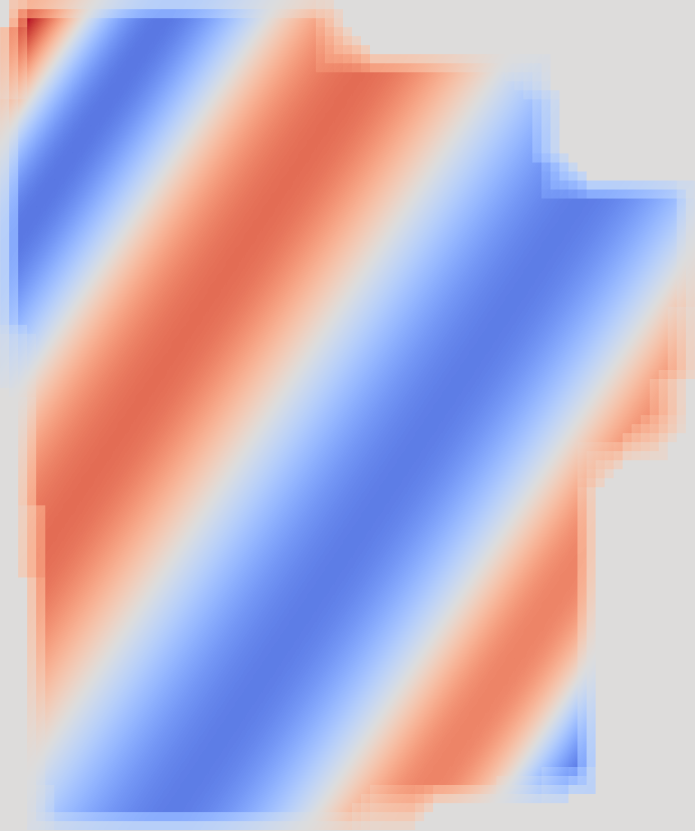}};%
  \end{tikzpicture}%
  \begin{tikzpicture}[baseline=(image.base)]%
  \node (image) at (0,0) {\includegraphics[width=.15\textwidth]{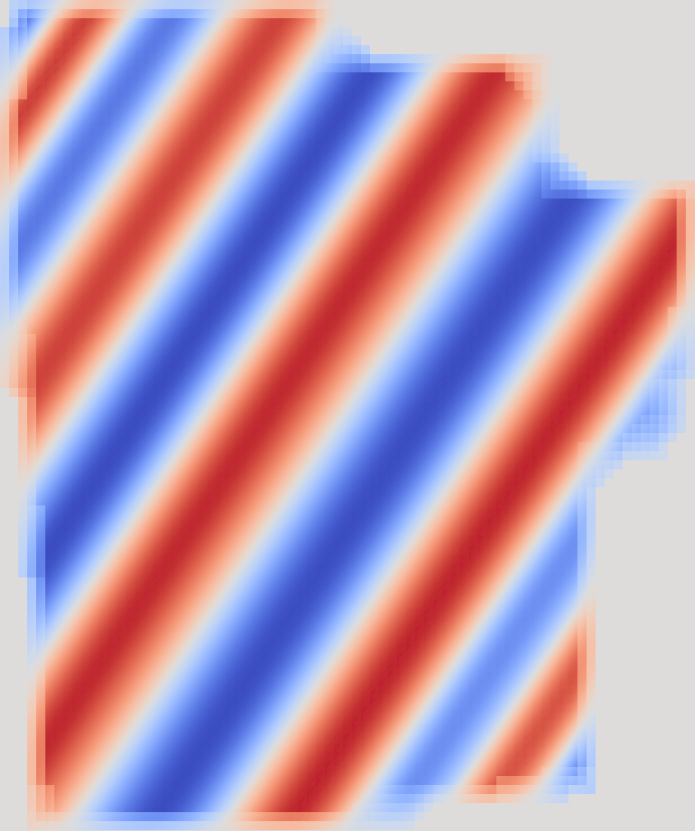}};%
  \end{tikzpicture}%
  \begin{tikzpicture}[baseline=(image.base)]%
  \node (image) at (0,0) {\includegraphics[width=.15\textwidth]{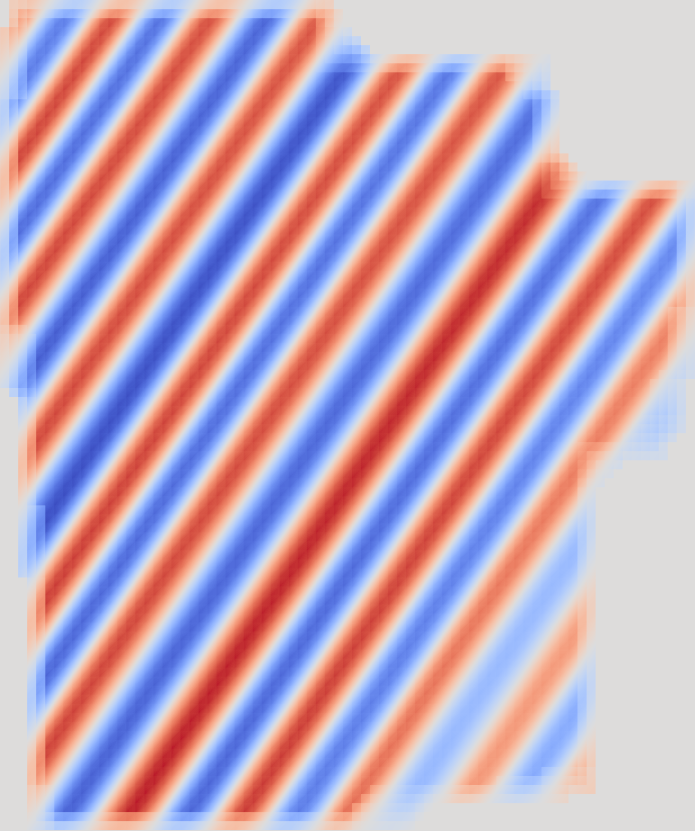}};%
  \end{tikzpicture}%
  \begin{tikzpicture}[baseline=(image.base)]%
  \node (image) at (0,0) {\includegraphics[width=.15\textwidth]{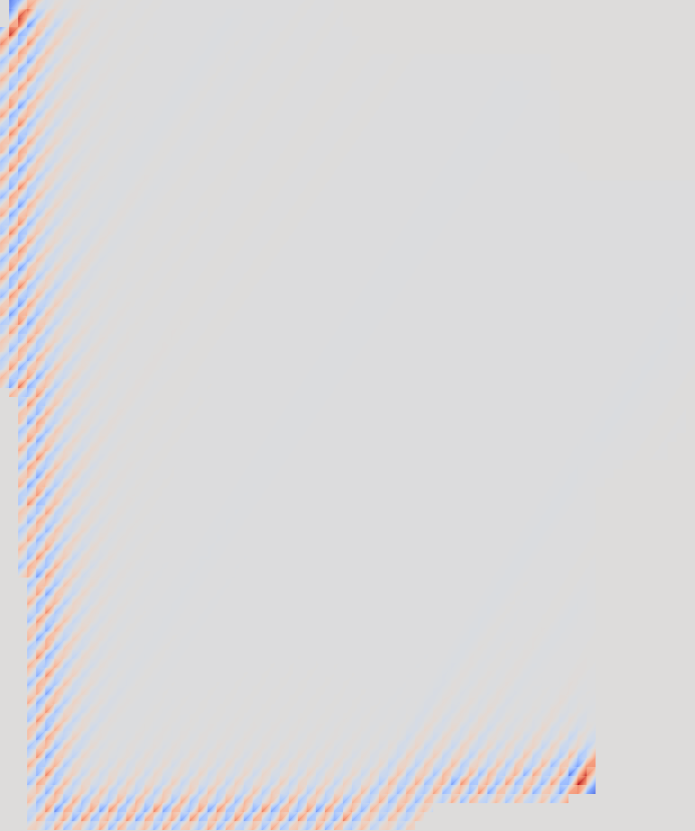}};%
  \end{tikzpicture}%
\end{center}
\caption{Selected eigenfunctions $\evecloc$ in an interior subdomain $\ovlpdomain_\ddindex$ of the homogeneous example from \Cref{sec:transport}. The model parameters are $\diffTensor\equiv10^{-3}$, $\velField\equiv0$ ($\peclet_h=0$) in the first row, and $\velField$ from \cref{eq:checkerboard:velField} with $\diffTensor\equiv10^{-3}$ ($\peclet_h\approx2$) and $\diffTensor\equiv10^{-6}$ ($\peclet_h\approx2\,000$) in the second and third row, respectively.}
\label{fig:transport:eigenfunctions:interior}
\end{figure}

To further investigate the influence of the grid P\'{e}clet number on the eigenfunctions and their decay, we visualize selected eigenfunctions on an interior subdomain $\ovlpdomain_\ddindex$ for three different examples: (a) pure diffusion, (b) diffusion with moderate convection, and (c) convection-dominated diffusion (\cref{fig:transport:eigenfunctions:interior}). In all cases, the constant function (corresponding to $\evalloc[1]=+\infty$) is found. While the eigenfunctions mainly vary around the boundary for the purely diffusive example, the influence and direction of the velocity field can clearly be seen in the other two rows. However, even in the examples including convection, the impact on the interior of the subdomains decreases the smaller the eigenvalues become. Comparing the cases with moderate and high grid P\'{e}clet numbers, it is apparent that this decay is delayed in the convection-dominated case, suggesting that more eigenfunctions are important for building the coarse space. This is both consistent with the coarse space dimensions we observe in \Cref{fig:transport} and the constants occurring in the eigenvalue decay from \cref{th:msgfem:evaldecay}. Nevertheless, as seen before, we obtain a robust preconditioner with moderately sized coarse spaces even in the convection-dominated regime.

\section{Conclusions}
\label{sec:conclusions}

By rewriting MS-GFEM as an iterative method, we have formulated a two-level hybrid restricted additive Schwarz preconditioner for convection-diffusion problems at high P\'{e}clet numbers. 
In extensive numerical experiments with highly heterogeneous diffusion on up to $10^5$ subdomains, we have demonstrated the preconditioner's robustness with respect to changes in the grid P\'{e}clet number in terms of GMRES iteration counts and coarse space dimensions. Since the method itself does not rely on strong assumptions on the specific discretization, we have obtained a flexible preconditioner that is applicable for different fine-scale discretizations on structured or unstructured meshes. Furthermore, we have observed robustness in a numerical experiment even in the vanishing-diffusion limit, which is not covered by the theory.

While we have seen in 3D simulations on a complex domain that the MS-GFEM preconditioner is more robust than other solvers, the computational cost can be high. In particular, solving the generalized eigenproblems in the setup phase clearly dominates the runtime. A promising approach to alleviating this is to follow \cite{Alber2025} and formulate a similar method with eigenproblems on rings instead of entire oversampling domains. Additionally, it is an interesting direction to investigate the influence of inexact iterative eigensolvers on the performance of the preconditioner. Since we do not use MS-GFEM as a multiscale method but as a preconditioner within GMRES, it is possible that the method is more forgiving of inaccurate eigenfunctions and a further speedup can be gained by exploiting this. Lastly, we aim to incorporate these ideas into a high-performance implementation of MS-GFEM-based preconditioners that are applicable to a wider range of problem types.


\appendix

\section*{Acknowledgments}
This work is funded by the Deutsche Forschungsgemeinschaft (DFG, German Research Foundation) under Germany's Excellence Strategy EXC 2181/1 - 390900948 (the Heidelberg STRUCTURES Excellence Cluster).
The authors thank Christian Alber and Chupeng Ma for fruitful discussions.

\section*{Data availability}
Data will be made available on reasonable request.

\renewcommand*{\bibfont}{\normalfont\small}
\printbibliography
\end{document}